\documentclass[11pt]{article}

\usepackage{amsfonts,amssymb,amsmath,amsthm,epsfig,euscript}
\usepackage{graphicx}
\usepackage{epstopdf}
\newtheorem{theorem}{Theorem}
\newtheorem{lemma}[theorem]{Lemma}

\long\def\symbolfootnote[#1]#2{\begingroup
\def\thefootnote{\fnsymbol{footnote}}\footnote[#1]{#2}\endgroup}

\newcommand{\fig}[2]{\begin{figure}[ht]
\centering\includegraphics[scale=.66]{#1}
\caption{#2}
\label{figure:#1}
\end{figure}}

\title{Ranking and unranking trees with a given number or 
a given set of leaves}

\author{Jeffery B.  Remmel\\
\small Department of Mathematics\\[-0.8ex]
\small U.C.S.D., La Jolla, CA, 92093-0112\\[-0.8ex]
\small \texttt{jremmel@ucsd.edu}\\
\and
S. Gill Williamson\\
\small Department of Computer Science and Engineering\\[-0.8ex]
\small U.C.S.D., La Jolla, CA, 92093-0404\\[-0.8ex]
\small \texttt{gwilliamson@ucsd.edu}
}

\date{ 
\small MR Subject Classifications: 05A15, 05C05, 05C20, 05C30}

\begin{document}
\maketitle

\begin{abstract}
In this paper, we provide algorithms to rank and unrank
certain degree-restricted classes of Cayley trees.  Specifically, we
consider classes of
trees that have a given set of leaves or that have a fixed number $k$ of 
leaves. Using properties of a bijection due to E\u gecio\u glu and
Remmel \cite{ER1}, we  reduce the problem of ranking and unranking these classes of degree-restricted trees to corresponding problems of 
ranking and unranking certain classes of set partitions. For fixed $k$, the number of Cayley trees with $n$ vertices and $k$ leaves grows roughly as $n!$ and hence 
the ranks have $O(nlog_2(n))$ bits. Our ranking and unranking algorithms require at most $O(n^2)$ comparisons of 
numbers $y \leq n$ plus $O(n)$ operations of multiplication, division, addition, substraction and comparision on numbers 
$x$ of length $O(nlog(n))$.
\end{abstract}

\section{Introduction}
In computational combinatorics, it is important to be able to {\it
efficiently\/} rank,
unrank, and randomly generate (uniformly) basic classes of combinatorial
objects. A ranking
algorithm for a finite set $S$ is a bijection from $S$ to the set 
$\{0, \cdots,|S|-1\}$.  An unranking algorithm is the inverse of a 
ranking algorithm.
Ranking and unranking techniques are useful for storage and 
retrieval of elements of
$S$.  Uniform
random generation plays a role in Monte Carlo methods and in search
algorithms such as hill
climbing or genetic algorithms over classes of combinatorial objects.
Uniform random generation of objects is always possible if 
one has an unranking algorithm since one can generate, 
uniformly, an integer in $\{0, \cdots, |S|-1\}$ and unrank.

We consider the set
$C_n$ of trees with vertex set $[n] = \{1, \ldots, n\}$.  These trees are sometimes called
Cayley trees and can
be viewed as the set of spanning trees of the complete graph $K_n$.
Ranking and unranking
algorithms for the set $C_n$ have been described by many authors. Indeed,
efficient ranking
and unranking algorithms have been given for classes of trees and forests
that considerably
generalize the Cayley trees (e.g., [3], [4], [5], [6], [7]).

In a previous paper \cite{RW1}, we considered a more refined problem, namely, the problem of ranking and unranking subsets of $C_n$ with a specified degree sequences or a specified multiset of degrees.  
Let $\vec{C}_{n,1}$ be the set of directed trees on 
$V$ that are rooted at 1. That 
is, a directed tree $T \in \vec{C}_{n,1}$ has all its edges 
directed towards its root 1. We replace $C_n$ with the equivalent set
$\vec{C}_{n,1}$.  For any tree $T \in C_n$,
$\sum_{i=1}^n deg_T(i) = 2n-2$.
If $\vec{s} = \langle s_1, \ldots, s_n \rangle$ is a sequence of
positive integers such
that $\sum_{i=1}^n s_i =2n-2$, then we  let 
$\vec{C}_{n,\vec{s}} = \{T \in 
\vec{C}_{n,1}:
\langle deg_T(1),
\ldots, deg_T(n)\rangle = \vec{s}\}$. Remmel and Williamson \cite{RW1} proved 
that  
\begin{equation}\label{count1}
|\vec{C}_{n,\vec{s}}| = \binom{n-2}{{s_1-1, \ldots,s_n-1}}.
\end{equation}
Similarly if $S =\{1^{\alpha_1}, \ldots,
(n-1)^{\alpha_{n-1}}\}$ is a multiset such that $\sum_{i=1}^{n-1} \alpha_i
\cdot i = 2n-2$ and $\sum_{i=1}^n \alpha_i = n$,
then we let
$\vec{C}_{n,S} = \{T \in \vec{C}_{n,1}: \{deg_T(1), \ldots, deg_T(n)\} = S \}$.  It is easy to see from (\ref{count1}) that 
\begin{equation}\label{count2}
|C_{n,S}| = \binom{n}{{\alpha_1, \ldots, \alpha_n}}\binom{n-2}{{s_1-1, \ldots,s_n-1}}.
\end{equation}
 
The basis of the ranking and unranking algorithms in \cite{RW1} 
for $\vec{C}_{n,\vec{s}}$ or $\vec{C}_{n,S}$ 
hinged on certian special properties of a 
bijection $\Theta$ between $\vec{C}_{n,1}$ and 
the class of functions ${\cal F}_n = 
\{f:\{2, \ldots, n-1\} \rightarrow [n]\}$ defined by  E\u gecio\u glu and Remmel \cite{ER1}. That is, in \cite{RW1}, we proved that  
for any vertex $i$, $1 +|f^{-1}(i)|$ equals 
that degree of $i$ in the tree $T= \Theta(f)$ when $\Theta(f)$ is considered 
as an undirected graph. This property allowed us to reduce the problem 
or ranking and unranking trees in $\vec{C}_{n,\vec{s}}$ or $\vec{C}_{n,S}$ to 
the problem of ranking and unranking certain classes of set partitions of 
$[n]$.  We were then able to modify known techniques for 
ranking and unranking set partitions \cite{W,NW} to construct efficient 
ranking and unranking algorithms for $\vec{C}_{n,\vec{s}}$ or $\vec{C}_{n,S}$.

In this paper, we shall give efficient algorithms to rank and unrank 
two other natural subsets of $\vec{C}_{n,1}$, namely, the set of trees which have a given number of leaves or a prespecified set of leaves.  
If $G = (V,E)$ is digraph and $v \in V$, we let
\begin{description}  
\item $indeg_G(v) =|\{u:(u,v) \in E\}|$, 
\item $outdeg_G(v) =|\{u:(v,u) \in E\}|$ and 
\item $deg_G(v) = indeg_G(v) + outdeg_G(v)$.
\end{description} 
If $T \in \vec{C}_{n,1}$, then we 
say that $i$ is a {\em leaf} of $T$ if and only if $deg_T(i) =1$. 
Fix $k$ such that $2\leq k \leq n-1$ 
and $\vec{C}_{n,1}^k$ equal the set of trees $T$ in $\vec{C}_{n,1}$ with $k$ leaves. Similarly, 
if $L$ is any subset of $\{1,\ldots,n\}$ of size $k$, we let 
$\vec{C}_{n,1}^{k,L}$ equal the set of trees $T \in \vec{C}_{n,1}$ 
such that $i$ is a leaf of $T$ if and only if $i \in L$.  
The main goal of this paper is to construct efficient 
ranking and unranking algorithms for the sets 
$\vec{C}_{n,1}^{k,L}$ or $\vec{C}_{n,1}^k$.

Just as in the case of the construction of the ranking and unranking algorithms 
for $\vec{C}_{n,\vec{s}}$ or $\vec{C}_{n,S}$ \cite{RW1}, 
the E\u gecio\u glu and Remmel bijection $\Theta$ 
allows us to reduce the problem 
problem ranking and unranking algorithms $\vec{C}_{n,1}^{k,L}$ or $\vec{C}_{n,1}^k$ 
to the problem of finding ranking and unranking certain classes of set 
partitions. That is,  we can restate the fundamental 
property of the bijection $\Theta: {\cal F}_n \rightarrow \vec{C}_{n,1}$ as 
\begin{equation}\label{fund}
deg_{\Theta(f)}(i) = |f^{-1}(i)|+1
\end{equation}
for all $f \in {\cal F}_n$ and $i \in [n]$. 
Now suppose that $L$ is a subset of $[n]$ of size $k$ where $2 \leq k \leq n$. 
Then by (\ref{fund}), it follows that if $T \in \vec{C}_{n,1}^{k,L}$ and 
$f = \Theta^{-1}(T)$, then $f^{-1}(i) = \emptyset$ if and only if $i \in L$. 
Thus if $J = [n]-L = \{j_1 < \ldots < j_{n-k}\}$, then 
$\langle f^{-1}(j_1), \ldots, f^{-1}(j_{n-k})\rangle$ must be an ordered set 
partition of $\{2, \ldots, n-1\}$ into $n-k$ nonempty parts. Conversely, 
if we are given an ordered set partition 
$\pi = \langle \pi_1, \ldots, \pi_{n-k}\rangle$ of $\{2, \ldots, n-1\}$ into $n-k$ non-empty parts, 
we can define a function $f \in {\cal F}_n$ such that 
$f^{-1}(i) = \emptyset$ if and only if $i \in L$ by setting 
$f^{-1}(j_t) = \pi_t$ for $t =1, \ldots, n-k$.  It   
follows that $|\vec{C}_{n,1}^{k,L}|$ is equal to the number of ordered set 
partitions of $\{2, \ldots, n-1\}$ into $n-k$ parts. We then develop 
algorithms for ranking and unranking such classes of ordered 
set partitions by modifying  ranking and unranking algorithms for the decreasing functions, 
permutations, and unordered set partitions found in \cite{W}.

Let $S_{n,k}$ denote the number of 
unordered set partitions of $[n]$ into $k$ parts.  The numbers 
$S_{n,k}$ are call the Stirling numbers of the second kind and 
they satisfy the following recursion:  
\begin{eqnarray*}
S_{n,k} &=& 0 \ \mbox{if either $k > n$ or $n <0$}, \\
S_{0,0} &=& 1, \ \mbox{and} \\
S_{n+1,k} &=& S_{n,k-1} + k S_{n,k}.
\end{eqnarray*}
Table 1 below gives the values of $S_{n,k}$ for $1 \leq n \leq 9$. 

\begin{center}
{\normalsize 
\begin{tabular}{|l||l|l|l|l|l|l|l|l|l|}
\hline 
 $\begin{matrix} \ & m \\ n & \ \end{matrix}$& 1 \ \ \ & 2 \ \ \ & 3 \ \ \  &4 \ \ \  & 5 \ \ \ & 6 \ \ \ & 7 \ \ \ & 8\ \ \ & 1\ \ \ \\
\hline \hline
1&1&&&&&&&&\\
\hline 
2&1&1&&&&&&& \\
\hline
3&1&3&1&&&&&&\\
\hline
4&1&7&6&&&&&&\\
\hline
5&1&15&25&10&1&&&&\\
\hline
6&1&31&90&65&15&1&&&\\
\hline 
7&1&63&301&350&140&21&1&&\\
\hline
8&1&127&966&1701&1050&266&28&1&\\
\hline
9&1&255&3025&7770&6951&2646&462&36&1\\
\hline
\end{tabular}}
\end{center}
\begin{center}
{\bf Table 1} The values of $S_{n,k}$
\end{center}

It follows from our arguments above that  
\begin{equation}\label{count21}
| \vec{C}_{n,1}^{k,L}| = (n-k)!S_{n-2,n-k}.
\end{equation}
Similarly for any $k$ such that $2\leq k \leq n-1$, it is easy to see 
that $\vec{C}_{n,1}^k$ is the disjoint union of all $C_n^{k,L}$ such that 
$L \subseteq [n]$ of size $k$ and hence   
\begin{equation}\label{count22}
| \vec{C}_n^{k}| = \binom{n}{k}  (n-k)!S_{n-2,n-k}.
\end{equation}
Note that both $|\vec{C}_{n}^{k,L}|$ and $|\vec{C}_{n}^k|$ can be as large as 
as $O(n!)$ so that the numbers involved in ranking and unranking can  
require $O(nlog(n))$ bits. 
We show that our ranking and unranking algorithms
 require at most $O(n^2)$ comparisons of 
numbers $y \leq n$ plus $O(n)$ operations of multiplication, division, addition, substraction and comparision on numbers 
$x < |\vec{C}_{n}^{k,L}|$ ($x < |\vec{C}_{n}^k|$).

The outline of this paper is as follows. In Section~2, we describe
the bijection $\Theta:{\cal F}_n \rightarrow \vec{C}_{n,1}$ of \cite{ER1}
and discuss some
of its key properties.  In Section~3, we show that both $\Theta$ and
$\Theta^{-1}$
can be computed in linear time. This result allows us to reduce the problem of
efficiently ranking and unranking trees in 
$\vec{C}_{n,\vec{s}}$ or $\vec{C}_{n,S}$ 
to the
problem of efficiently ranking and unranking certain classes of ordered set
partitions.  In section 4, we shall recall the algorithms due to Williamson \cite{W} for ranking and unranking decreasing functions, permutations and 
unordered set partitions which will be the building blocks of 
our final ranking and unranking algorithms.  
Finally in Section~5,  we  shall give our 
 ranking and unranking algorithms for 
the sets $\vec{C}_{n}^{k,L}$ or $\vec{C}_{n}^k$ and 
give examples.

\section{The $\Theta$ Bijection and its Properties}

In this section, we shall review the bijection 
$\Theta: {\cal F}_n \rightarrow \vec{C}_{n,1}$ due to 
E\u gecio\u glu and Remmel \cite{ER1} and give some of its properties. 

Let $[n] = \{1, 2, \ldots, n\}$. For each function 
$f:\{2, \ldots, n-1\} \rightarrow [n]$, we associate a directed graph $f$, 
$graph(f) = ([n],E)$ by setting $E = \{\langle i,f(i) \rangle :i = 2, \ldots, n-1\}$.   Following \cite{RW}, given  any  directed edge $(i,j)$
where $1 \leq i,j \leq n$, we define the weight of $(i,j)$, $W((i,j))$, by
\begin{equation}\label{wedge}
W((i,j)) =
\left\{ \begin{array}{ll}
p_is_j \  \mbox{if $i < j$}, \\
q_it_j \  \mbox{if $i \geq j$}
\end{array} \right.
\end{equation}
where $p_i, q_i, s_i, t_i$ are variables for $i = 1, \ldots, n$.
We shall call a directed edge $(i,j)$ a {\em descent edge} if $i \geq j$ and an
{\em ascent edge} if $i < j$.  We then define the weight of any digraph $G =
([n], E)$ by
\begin{equation}\label{wgraph}
W(G) = \prod_{(i,j) \in E} W((i,j)).
\end{equation}
A moment's thought will convince one that, in
general, the digraph corresponding to a function $f \in {\cal F}_n$
will consists of $2$ root-directed trees
rooted at vertices $1$ and $n$ respectively, with all
edges directed toward their roots, plus a number of
directed cycles of length $\geq 1$. For each vertex $v$ on a given cycle,
there is possibly a root-directed tree
attached to $v$  with $v$ as the root and all edges directed toward $v$.
Note the fact that there are trees rooted at vertices
$1$ and $n$ is due to the fact that these elements are not in the domain
of $f$.  Thus 
there can be no directed edges out of any of these vertices. We let the
weight of
$f$, $W(f)$, be the weight of the digraph $graph(f)$ associated with $f$.

To define the bijection $\Theta$, we first imagine
that the directed graph corresponding
to $f \in {\cal F}$ is drawn  so that
\begin{description}
\item[(a)] the trees rooted at $n$ and $1$ are drawn on the extreme left
and the extreme
right respectively with their edges directed upwards,

\item[(b)] the cycles are drawn so that their vertices form a directed path
on the line
between $n$ and $1$, with one back edge above the line, and  the
root-directed tree attached
to any vertex on a cycle is
drawn below the line between $n$ and $1$  with its edges directed upwards,

\item[(c)] each cycle $c_i$ is arranged so that its maximum
element $m_i$ is on the right, and 

\item[(d)] the cycles are arranged from left to right by decreasing maximal elements.
\end{description}
Figure~1 pictures a function $f$ drawn according to the rules  
(a)-(d) where $n =23$.

\fig{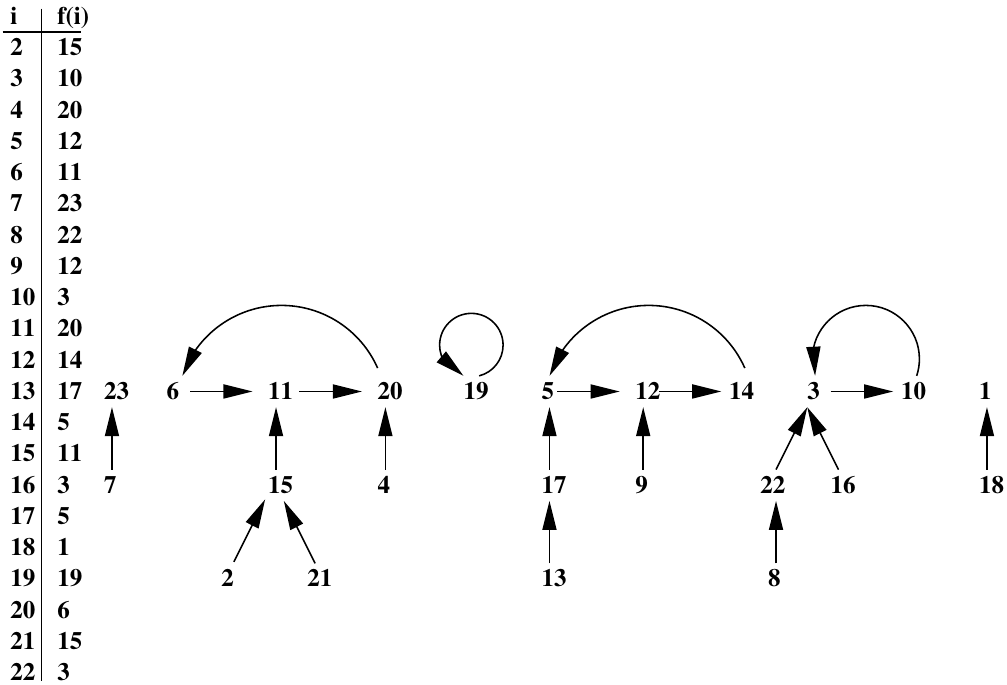}{The digraph of a function.}


This given, suppose that the digraph of $f$ is drawn as described above and
the cycles of $f$ are $c_1(f), \ldots, c_a(f)$, reading from left to right.
We let $r_{c_i(f)}$ and $l_{c_i(f)}$ denote the right and left endpoints of
the cycle $c_i(f)$ for $i= 1, \ldots,a$. Note that if $c_i(f)$ is a
1-cycle, then
we let $r_{c_i(f)} = l_{c_i(f)}$ be the element in the 1-cycle.
$\Theta(f)$ is obtained from $f$ by simply deleting the back edges
$(r_{c_i(f)},l_{c_i(f)})$ for $i = 1, \ldots, a$ and adding the directed
edges $(r_{c_i(f)},l_{c_{i+1}(f)})$ for $i =
1, \ldots, a-1$ plus the directed edges $(n,l_{c_1(f)})$ and
$(r_{c_a(f)},1)$.
That is, we remove all the back
edges that are above the line, and then we connect $n$ to the lefthand
endpoint of the first cycle, the righthand
endpoint of each cycle to the lefthand endpoint of the cycle following it,
and we connect the righthand endpoint of
the last cycle to $1$.
For example, $\Theta(f)$ is pictured in Figure~2  for the $f$ given in
Figure~1.
If there are no cycles in $f$, then $\Theta(f)$ is simply the result of
adding the directed edge $(n,1)$ to the digraph of $f$.

\fig{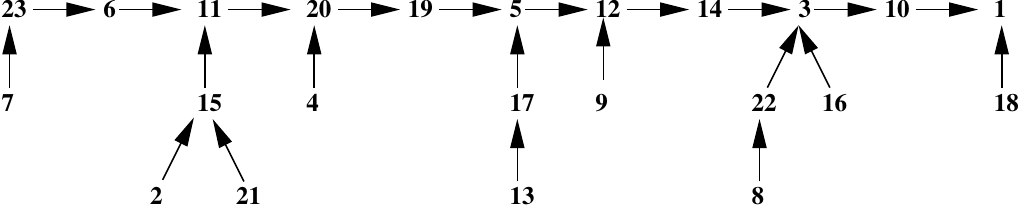}{$\Theta(f)$.}


To see that $\Theta$ is a bijection, we shall describe how to define
$\Theta^{-1}$. The key observation is that we need only recover 
that the directed edges $(r_{c_i(f)},l_{c_{i+1}(f)})$ for $i =
1, \ldots, a-1$.  However it is easy to see that 
$r_{c_1(f)} = m_1$ is the largest element on the path from $n$ to $1$ 
in the tree $\Theta(f)$.  That is, $m_1$ is then largest element in its cycle and by definition, it is larger than all the largest elements in any other cycle so 
that $m_1$ must be the largest interior element on the path from $n$ to 1. 
Then by the same reasoning, $r_{c_2(f)} = m_2$ is the largest element 
on the path from $m_1$ to 1, etc. Thus we can find $m_1, \ldots, m_t$.    
More formally, given a tree  $T\in \vec{C}_{n,1}$, consider the path
$$m_0=n, x_1, \ldots, m_1, x_2, \ldots m_2, \ldots, x_t, \ldots , m_t, 1$$
where  $m_i$ is the maximum interior vertex on the path
from $m_{i-1}$ to $1$, $1\leq i\leq t$.  If $(m_{i-1},m_i)$ is an edge on
this path, then it is
understood that $x_i, \ldots , m_i = m_i$ consists of just one vertex and
we define $x_i=m_i$.  Note that by definition $m_0 = n > m_1 > \ldots > m_t$.
We obtain the digraph $\Theta^{-1}(T)$ from $T$ via the following
procedure.\\
\ \\
{\bf Procedure for computing $\Theta^{-1}(T):$}

\medskip\noindent
(1) First we declare that any  edge $e$ of $T$
which is not an edge of the path
from $n$ to $1$ is
an edge of $\Theta^{-1}(T)$.\\
\ \\
(2) Next we remove all edges of the form
$(m_t, 1)$  or $(m_{i-1}, x_i)$ for $1\leq i\leq t$.\\
\ \\
Finally for each $i$ with $1\leq i\leq t$, we consider the subpath
$x_i, \ldots , m_i$.\\
\ \\
(3) If $m_i=x_i$, create a directed loop
$(m_i, m_i)$.\\
\ \\
(4) If $m_i\neq x_i$,
convert the
subpath
$x_i, \ldots , m_i$ into the \\
directed cycle $x_i, \ldots, m_i, x_i$.\\
\ \\

Next we consider two important properties of the bijection $\Theta$.
First $\Theta$ has an important weight preserving property. 
We claim that if $\Theta(f) = T$, then 
\begin{equation}\label{weightpres}
q_nt_1 W(f) = W(T).
\end{equation}
That is, by our conventions, any backedge 
$(r_{c_{i}(f)},l_{c_{i}(f)})$ are descent edges so that its weight is
$q_{r_{c_{i}(f)}}t_{l_{c_{i}(f)}}$.  Thus the total weight of the backedges is
\begin{equation}\label{fnw}
\prod_{i=1}^a q_{r_{c_{i}(f)}}t_{l_{c_{i}(f)}}.
\end{equation}
Our argument above shows that all the new edges that we add are also
descent edges so that the weight of the new edges is
\begin{equation}
q_nt_{l_{c_1(f)}} (\prod_{i=1}^{a-1} q_{r_{c_{i}(f)}}t_{l_{c_{i+1}(f)}})
q_{r_{c_a(f)}}t_1 = q_nt_1 \prod_{i=1}^a q_{r_{c_{i}(f)}}t_{l_{c_{i}(f)}}.
\end{equation}
Since all the remaining edges have the same weight in both the digraph of
$f$ and in the digraph $\Theta(f)$, it follows that $q_nt_1W(f) = W(\Theta(f))$ as claimed.

It is easy to see that 
\begin{equation}
\sum_{f \in {\cal F}_n} W(f) = \prod_{i=2}^{n-1} [q_i(t_1 + \cdots + t_i) + p_i(s_{i+1} + \cdots + s_n)].
\end{equation}
Thus we have the following result which is implicit in \cite{ER1} and it explicit in \cite{RW}.
\begin{theorem}\label{weightthm1}
\begin{equation}
\sum_{T \in \vec{C}_{n,1}} W(T) = q_n t_1 \prod_{i=2}^{n-1} 
[q_i(t_1 + \cdots + t_i) + p_i(s_{i+1} + \cdots + s_n)].
\end{equation}
\end{theorem}

Next we turn to a second key property of the $\Theta$ bijection. It 
is easy to see from Figures 1 and 2 that 
deleting the back edges
$(r_{c_i(f)},l_{c_i(f)})$ for $i = 1, \ldots, a$ in $graph(f)$ and adding the directed edges $(r_{c_i(f)},l_{c_{i+1}(f)})$ for $i =
1, \ldots, a-1$ plus the directed edges $(n,l_{c_1(f)})$ and
$(r_{c_a(f)},1)$ to get $\Theta(f)$ does not change the indegree of 
any vertex except vertex 1. 
That is,    
\begin{equation}\label{Deg1}
indeg_{graph(f)}(i) = indeg_{\Theta(f)}(i)\ \mbox{for} \ i =2, \ldots ,n.
\end{equation}
It is also easy to see that in going from $graph(f)$ to $\Theta(f)$, the 
indegree of vertex 1 increases by 1, i.e., 
\begin{equation}\label{Deg2}
1+ indeg_{graph(f)}(i) = indeg_{\Theta(f)}(1).
\end{equation}
When we consider $\Theta(f)$ as an undirected graph $T$, then it is 
easy to see that $deg_T(i) = outdeg_{\Theta(f)} + indeg_{\Theta(f)}$. Thus 
since the outdegree of $i$ in $\Theta(f)$ is 1 if $ i \neq 1$ and the outdegree of 1 in $\Theta(f)$ is zero, equations (\ref{Deg1}) and (\ref{Deg2}) imply 
the following theorem.

\begin{theorem}\label{degthm}
 Suppose that $T$ is the undirected tree corresponding 
to $\Theta(f)$ where $f \in {\cal F}_n$, then for $i =1, \ldots, n$, 
\begin{equation}\label{Deg3}
deg_T(i) = 1 + |f^{-1}(i)|.
\end{equation}
\end{theorem}
{\em Proof} By our definition of $graph(f)$, it follows that 
$indeg_{graph(f)}(i) = |f^{-1}(i)|$ for $i =1, \ldots, n$. 
Thus by (\ref{Deg1}), for $i = 2, \ldots, n$,
\begin{eqnarray*}
deg_T(i) &=&  outdeg_{\Theta(f)}(i) + indeg_{\Theta(f)}(i) \\
&=& 1 + indeg_{\Theta(f)}(i) \\
&=& 1 + indeg_{graph(f)}(i) \\
&=& 1 + |f^{-1}(i)|.
\end{eqnarray*}
Similarly by (\ref{Deg2}), 
\begin{eqnarray*}
deg_T(1) &=&  outdeg_{\Theta(f)}(1) + indeg_{\Theta(f)}(1) \\
&=& 0 + indeg_{\Theta(f)}(1) \\
&=& 1 + indeg_{graph(f)}(1) \\
&=& 1 + |f^{-1}(1)|.
\end{eqnarray*}
\hfill $\Box$

\section{Construction of the spanning forests from the the function
table in time $O(n)$}

In this section, we shall briefly outline the proof that one can 
compute the bijections  $\Theta$ and its inverse in linear time.  
Suppose we are given $f \in {\cal F}_n$. 
Our basic data structure for the  function $f$ is a list of pairs $
\langle i,f(i) \rangle$ for $i = 2, \ldots, n-1$.   
Our goal is to construct the directed 
graph of $\Theta(f)$ from our data structure for $f$, that is, 
for $i = 1, \ldots, n$, we want to find the set of pairs, $\langle i,t_i \rangle$,  such that there is directed edge from $i$ to $t_i$ in $\Theta(f)$. We shall prove the following. 

\begin{theorem}\label{linbij}
We can compute the bijection
$\Theta:{\cal F}_n \rightarrow \vec{C}_{n,1}$ 
and its inverse in linear time.
\end{theorem}
{\em Proof.}
We shall not try to give the most efficient algorithm to construct 
$\Theta(f)$ from $f$. Instead, we shall give an 
outline the basic procedure which shows that one can 
construct $\Theta(f)$ from $f$ in linear time. For 
ease of presentation, we shall organize our procedure so that it 
makes four linear time passes through the basic data structure for $f$ to 
produce the data structure for $\Theta(f)$.\\ 
\ \\
{\bf Pass 1.} {\it Goal: Find, in linear time in $n$, a set of representatives
$t_1, \ldots, t_r$ of the cycles of the directed graph of the function $f$.} \\
To help us find $t_1, \ldots, t_r$, we shall maintain an array 
$A[2], A[3], \cdots A[n-1]$,  
where for each $i$, $A[i]=(c_i, p_i, q_i)$ is a triple of integers such 
$c_i \in \{0,\ldots, n-1\}$ and $\{p_i,q_i\} \subseteq 
\{-1,2,\cdots, n-1\}$.  
The $c_i$'s will help us keep track of what loop we are in relative to 
the sequence of operations described below. Then our idea is to  
maintain, through the $p_i$ and $q_i$, a doubly linked list of the
locations $i$ in $A$ where $c_i=0$, 
and we obtain pointers to the first and last elements of this
doubly linked list. It is a standard exercise that these 
data structures can be maintained in
linear time.

Initially, all the $c_i$'s will be 
zero.  In general, if $c_i=0$, then $p_i$ will be the largest integer $j$ such that  
$2 \leq j < i$ for which $c_{j}=0$ if there is such a $j$ and 
$p_i = -1$ otherwise. 
Similarly, we set 
$q_i>i$ to be the smallest integer $k$ such that $n-1 \geq k > i$  for which
$c_{k}=0$ if there is such a $k$ and 
$q_i=-1$ if there is no such $k$.   If $c_{2}>0$, then $q_{2}$ is
the smallest integer $j > 2$  such that $c_{j}=0$ and 
${q_{2}}=-1$ if there is  no such integer $j$.  If $c_{n-1}>0$,  then
${p_{n-1}}$ is the largest integer $k < n_1$ such that $c_{k}=0$ and
${p_{n-1}}=-1$ if there is  no such integer $k$.

We initialize $A$ by setting $A[2]=(0,-1,q_{2})$, 
$A[i]=(0,i-1,i+1)$ for $m+1<i<n-1$, and $A[n-1]=(0,p_{n-1},-1)$.
If $2<n-1$ then $q_{2}=3$ and $p_{n-1}=n-2$.  Otherwise ($2=n-1$),
these quantities are both $-1$.\\
\ \\
{\bf LOOP(1):} Start with $i_1=2$, setting $c_{2}=1$. 
Compute   $f^0(2), f^1(2),  f^2(2), \ldots, f^{k_1}(2)$, each
time updating
$A$ by setting $c_{f^{j}(2)}=1$ and adjusting pointers, until, 
prior to setting $c_{f^{k_1}(2)}=1$, we discover that  either \\
\ \\
(1) $f^{k_1}(2) \in \{1, n\}$, in which case  
we have reached a node in ${graph}(f)$ which is not in the domain of $f$ and 
we start over again with the
$2$ replaced by the smallest $i$ for which  $c_i=0$, or \\
\ \\
(2) $x=f^{k_1}(2)$ already satisfies $c_x=1$.  This condition indicates that
the value $x$ has
already occurred in the sequence $2, f(2),  f^2(2), \ldots,
f^{k_1}(2)$.  Then we set $t_1=f^{k_1}(2)$.\\
\ \\
{\bf LOOP(2):} Start with $i_2=q_{m+1}$ which is the location of the first $i$ such that $c_i=0$, and repeat the calculation of LOOP1 with
$i_2$ instead of $i_1=2$.
In this manner, generate $f^0(i_2), f^1(i_2),  f^2(i_2), \ldots,
f^{k_2}(i_2)$, each time updating
$A$ by setting $c_{f^{j}(i_2)}=2$ and adjusting pointers,
until either \\
\ \\
(1) $f^{k_1}(i_2) \in \{1,n\}$, in which case we have 
reached a node in
${graph}(f)$ which is not in the domain of $f$ and we start
over again with the
$i_2$ replaced by the smallest $i$ for which  $c_i=0$, or \\
\ \\
(2) $x=f^{k_1}(i_2)$ already satisfies $c_x=2$. (This condition indicates that
the value $x$ has
already occurred in the sequence $i_2, f(i_2),  f^2(i_2), \ldots,
f^{k_1}(i_2)$.) Then we set $t_2=f^{k_1}(i_2)$.\\
\ \\
We continue this process until $q_{2}=-1$. 
At this point, we will have generated
$t_1, \ldots, t_r$, where the last loop was LOOP($r$). 
The array $A$ will be such that, for all $2\leq i\leq n-1$, 
$1\leq c_i\leq r$ identifies the LOOP in which
that particular domain value $i$ occurred in our computation described above.\\
\ \\
{\bf Pass 2.} {\em Goal: For $i =1, \ldots, r$, find the largest element $m_i$ in the cycle determined by $t_i$.}\\ 
It is easy to see that this computation can be  done in linear time by
one pass through the array $A$ computed in Pass 1 above. At the end of 
Pass 2, we set $l_i = f(m_i)$. Thus when we draw the cycle containing $t_i$ according to our definition of $\Theta_j(f)$, $m_i$ will be right most element in the and $l_i$ will be the left most element of the cycle containing $t_i$. However, at this point, we have not ordered the cycles appropriately.  This ordering will be done in the next pass. \\
\ \\
{\bf Pass 3.} {\em Goal:  Sort $(l_1,m_1), \ldots, (l_k,m_k)$ so that they are
appropriately ordered according the criterion for the bijection $\Theta(f)$ 
as described in by condition (a) -(d)}.\\
\ \\
Since we order the cycles from left to right according to decreasing 
maximal elements, it is then easy to see that our desired 
ordering can be constructed via a lexicographic bucket sort. 
(See Williamson's book \cite{W} for details on the fact that a 
lexicographic bucket sort can be carried out in linear time.)\\
\ \\ 
{\bf Pass 4.} {\em Goal: Construct the digraph of $\Theta(f)$ from 
the digraph of $f$.}\\
\ \\
We modify the table for $f$ to produce the table for $\Theta(f)$ as follows.
Assume that $(l_1,m_1), \ldots, (l_k,m_k)$ is the sorted list
coming out of Pass 3.
Then we modify the table for $f$ so that we add entries for 
the directed edges $\langle n, l_1\rangle $ and $ \langle m_k, 1 \rangle$ and
modify entries of the pairs starting with $m_1, \ldots, m_k$ so that their corresponding second elements are $l_2, \ldots, l_k, 1$ respectively.
This can be done in linear time using our data structures.

Next, consider the problem of computing the inverse of $\Theta$.
Suppose that we are given the data structure of the tree  $T \in \vec{C}_1$, i.e. we are given a 
set of pairs , $\langle i,t_i\rangle$,  such that there is a directed edge from $i$ to $t_i$ in $T$.
Recall that the computation of $\Theta^{-1}(T)$ consists of two basic steps. \\
\ \\
{\bf Step 1.} Given a tree  $T \in \vec{C}_{n,1}$, consider the path
$$m_0=n, x_1, \ldots, m_1, x_2, \ldots m_2, \ldots, x_t, \ldots , m_t, 1$$
where  $m_i$ is the maximum interior vertex on the path
from $m_{i-1}$ to $1$, $1\leq i\leq t$.  If $(m_{i-1},m_i)$ is an edge on
this path, then it is
understood that $x_i, \ldots , m_i = m_i$ consists of just one vertex and
we define $x_i=m_i$.  Note that by definition $m_0 = n > m_1 > \ldots > m_t$.\\
\ \\
First it is easy to see that by making one pass through the data structure for $F$, 
we can construct the directed path $n \rightarrow a_1 \rightarrow \ldots 
\rightarrow a_r$ where $1 = a_r$. 
In fact, we can construct a doubly linked list $(n,a_1, \ldots, a_{r-1},1)$ with pointers 
to the first and last elements in linear time. If we traverse the list in reverse order, 
$(1, a_{r-1}, \ldots, a_1,n)$, then it easy to see that $m_t = a_{r-1}$, $m_{t_1}$ is the 
next element in the list $(a_{r-2}, \ldots ,a_1)$ which is greater than $m_t$ and, in 
general, having found $m_i = a_s$, then $m_{i-1}$ is the first element in the list 
$(a_{s-1}, \ldots, a_1)$ which is greater than $m_i$.  Thus it is not difficult to see 
that we can use our doubly linked list to produce the factorization 
$$m_0=n, x_1, \ldots, m_1, x_2, \ldots m_2, \ldots, x_t, \ldots , m_t, 1 $$  
in linear time. \\
\ \\ 
{\bf Step 2.}
We obtain the digraph $\Theta^{-1}(T)$ from $T$ via the following
procedure.\\
\ \\
{\bf Procedure for computing $\Theta_j^{-1}(F):$}

\medskip\noindent
(1) First we declare that any  edge $e$ of $T$
which is not an edge of the path
from $n$ to $1$ is
an edge of $\Theta^{-1}(T)$.\\
\ \\
(2) Next we remove all edges of the form
$(m_t, 1)$  or $(m_{i-1}, x_i)$ for $1\leq i\leq t$.\\
\ \\
Finally for each $i$ with $1\leq i\leq t$, we consider the subpath
$x_i, \ldots , m_i$.\\
\ \\
(3) If $m_i=x_i$, create a directed loop
$(m_i, m_i)$.\\
\ \\
(4) If $m_i\neq x_i$, then ,
convert the
subpath
$x_i, \ldots , m_i$ into the \\
directed cycle $x_i, \ldots, m_i, x_i$.\\
\ \\
Again it is easy to see that we can use the data structure for 
$T$, our doubly linked list, and our path factorization, 
$m_0=n, x_1, \ldots, m_1, x_2, \ldots m_2, \ldots, x_t, \ldots , m_t, j$ to construct the 
data structure for $graph(f)$  where $f = \Theta^{-1}(T)$ in linear time.  
\hfil $\Box$.

Given that we can carry out the bijection $\Theta$ and its inverses in 
linear time, it follows that in linear time, we can reduce the problem of constructing  ranking and unranking algorithms for $C_n$ to the problem 
of constructing ranking and unranking algorithms for the corresponding function class ${\cal F}_n$.

\section{Machinary for Ranking and Unranking Algorithms.}

Fix $k$ such that $2\leq k \leq n-1$. Recall that $\vec{C}_{n,1}^k$ equals the set of trees $T$ in $\vec{C}_{n,1}$ with $k$ leaves. Similarly, 
if $L$ is any subset of $\{1,\ldots,n\}$ of size $k$, we let 
$\vec{C}_{n,1}^{k,L}$ equal the set of trees $T \in \vec{C}_{n,1}$ 
such that $i$ is a leaf of $T$ if and only if $i \in L$.  The main 
goal of this section is to develop the basic machinary that is needed to 
give our final 
ranking and unranking algorithms for 
the sets $\vec{C}_{n,1}^{k,L}$ or $\vec{C}_{n,1}^k$.

Our ranking and unranking algorithms for $\vec{C}_{n,1}^k$ 
are based on six reductions.
\begin{enumerate}
\item  By the E\u gecio\u glu and Remmel bijection $\Theta$ of section 2, 
the problem of ranking and unranking trees in $\vec{C}_{n,1}^k$ can 
be reduced to the problem of ranking and unranking the set $F_{n,k}$ 
of functions 
$f:{2, \ldots, n-1} \rightarrow [n]$ such that $|\{i \in [n]: f^{-1}(i) = \emptyset\}| = k$.  

\item To specify the set of $i \in [n]$ such that $f^{-1}(i) = \emptyset$ for 
an $f \in F_{n,k}$, we 
specify a decreasing function $g_f:\{1, \ldots, k\} \rightarrow [n]$ whose 
range if $\{i:f^{-1}(i) =\emptyset\}$.  
Let ${\cal DF}_{n,k}$ denote the set of decreasing functions $f:\{1, \ldots, k\} \rightarrow \{1, \ldots, n\}$. Then given $g \in {\cal DF}_{n,k}$, we let $F_{n,g}$ equal the set of functions 
$f:{2, \ldots, n-1} \rightarrow [n]$ such that $\{i \in [n]: f^{-1}(i) = \emptyset\} = range(g)$. 

\item Given $g_f \in {\cal DF}_{n,k}$, we specify a function $f \in F_{n,g}$, 
by giving an ordered set partition 
$\pi_f$ of $[n-2]$ into $n-k$ parts. That is, 
suppose that $\pi = \langle \pi_1, \ldots, \pi_{n-k}\rangle$ is an 
ordered set partition of $[n-2]$ into $n-k$ parts and 
$[n]-range(g) = \{i_1 < \cdots < i_{n-k}\}$. Let $\pi^+ = \langle \pi_1^+, \ldots, \pi_{n-k}^+\rangle$ denote the set partition of $\{2, \ldots, n-1\}$ which results from $\pi$ by replacing each element $x\leq n-2$ by $x+1$. Then 
we can specify an $f \in F_{n,g}$ by declaring that 
$f^{-1}(i_j) = \pi^+_j$ for $j =1, \ldots, n-k$.

\item To specify an ordered set partition 
of $[n-2]$ into $n-k$ parts, we shall specify an unordered 
set partition $\Gamma_f = \langle \Gamma_1, \ldots, \Gamma_{n-k}\rangle$ 
of $[n-2]$ into $n-k$ parts where $min(\Gamma_1) < \cdots < min(\Gamma_k)$ and 
a permutation $\sigma_f = \sigma_1 \ldots \sigma_{n-k}$ in the symmetric group $S_{n-k}$.  That is, 
given $\Gamma_f$ and $\sigma_f$, we let $\Gamma_\sigma$ be the ordered 
set partition of $[n-2]$ into $k$ parts where 
$\Gamma_\sigma = \langle \Gamma_{\sigma_1}, \ldots, \Gamma_{\sigma_{n-k}}\rangle$.

\item We shall associate to each permutation $\sigma = \sigma_1 \ldots \sigma_{n} \in S_{n}$, a sequence,  $h_\sigma$, called the 
direct insertion sequence  associated 
with $\sigma$.  We defined $h_\sigma = (h(1), \ldots ,h(n))$ by recursion 
as follows. First if $\sigma \in S_1$, $h_\sigma = (0)$.   
If $n \geq 0$ and $\sigma = \sigma_1 \ldots, \sigma_n \in S_n$ where 
$\sigma_1 =j$, then let $\sigma^-$ be the permutation derived 
from $\sigma_2 \ldots \sigma_n$ where we replace each $i > j$ in 
the sequence $\sigma_2 \ldots \sigma_n$  by $i-1$.  For example, 
if $\sigma = 4 \ 1 \ 2 \ 5 \ 6 \ 3$, then $\sigma^- = 1 \ 2 \ 4 \ 5 \ 3$. 
This given, suppose $h_{\sigma^-} = (h^-(1), \ldots, h^-(n-1))$, then 
we define $h_\sigma = (j-1, h^-(1), \ldots, h^-(n-1))$. For example, 
$\sigma = 4 \ 1 \ 2 \ 5 \ 6 \ 3$, then $h_\sigma = (3,0,0,1,1,0)$.

\item We shall associate to each unordered set partition 
$\Gamma = \langle \Gamma_1, \ldots, \Gamma_{n-k}\rangle$ 
of $[n-2]$ into $n-k$ parts where 
$min(\Gamma_1) < \cdots < min(\Gamma_{n-k})$, 
a sequence $s_\Gamma = (s(1), \ldots, s(n-2))$ associated to its 
corresponding restricted growth function. That is,  
$s(i) =j-1$ if and only if $i \in \Gamma_j$. For example, 
suppose $\Gamma = \langle \{1,4\}, \{2,3,8\},\{5\},\{6,7\}\rangle$, 
then $s_\gamma = (0,1,1,0,2,3,3,1)$. 
\end{enumerate}

It follows that we can specify any $f \in F_{n,k}$ by 
a triple $\langle g, \sigma, \Gamma \rangle$ where 
$g \in \in {\cal DF}_{n,k}$, 
$\sigma \in S_{n-k}$ and $\Gamma$ is an unordered set partition of 
$[n-2]$ into $n-k$ parts. We can thus identity $f$ with a triple of 
sequences 
$$Seq(f) = \langle (g(1), \ldots, g(k)), 
(h(1), \ldots, h(n-k)),(s(1), \ldots, s(n-2)) \rangle$$
where $h_\sigma = (h(1), \ldots, h(n-k))$ and $s_\Gamma = 
(s(1), \ldots, s(n-2))$.  We can then order the 
functions $f \in F_{n,k}$ according to the lexicographic order of 
their associated sequences $Seq(f)$.

\fig{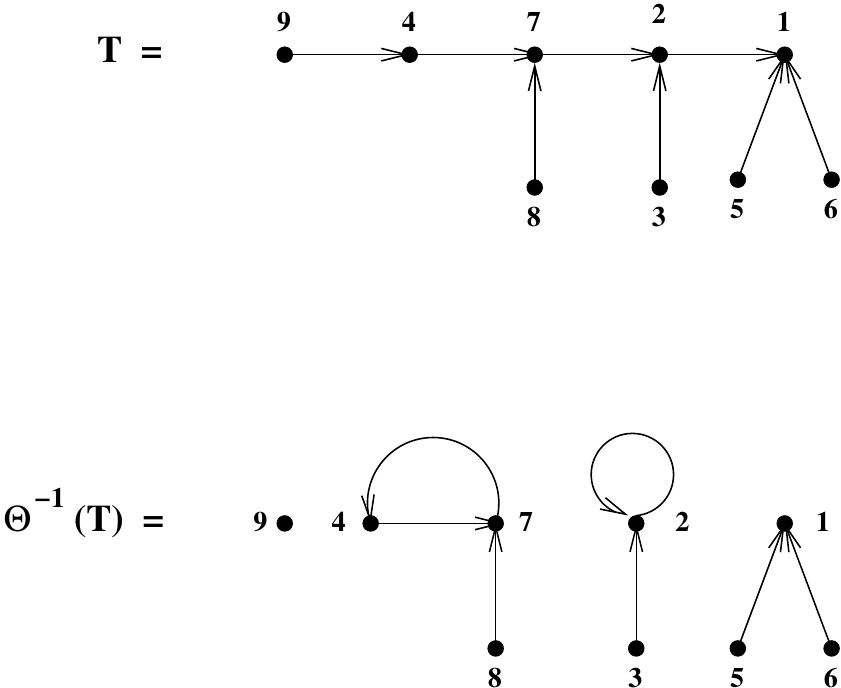}{A tree in $\vec{C}_{9,1}^5$.}


For example, consider 
the tree in $T \in \vec{C}_{9,1}^5$ pictured in Figure \ref{figure:Texample}. 
We have 
also pictured $\Theta^{-1}(T) =f$ so that $f$ can be specified 
by its table of preimages $f^{-1}(1), \ldots, f^{-1}(9)$ given in 
Table 2.\\
\ \\ 
{\normalsize 
\begin{tabular}{|l|l|l|l|l|l|l|l|l|l|}
\hline 
$f^{-1}(1)$ & $f^{-1}(2)$ & $f^{-1}(3)$ & $f^{-1}(4)$ & $f^{-1}(5)$ & 
$f^{-1}(6)$ & $f^{-1}(7)$ & $f^{-1}(8)$ & $f^{-1}(9)$ \\
\hline
$\{5,6\}$ & $\{2,3\}$ & $\emptyset$ & $\{7\}$ & $\emptyset$ & $\emptyset$ & 
$\{4,8\}$ & $\emptyset$ & $\emptyset$ \\
\hline 
\end{tabular}}
\begin{center}
{\bf Table 2}
\end{center}
Thus $f$ is associated with the decreasing functions $g_f \in {\cal DF}_{9,5}$ given 
by $g_f =(9,8,6,5,3)$ and the ordered set partition 
$\pi^+_f = \langle \{5,6\}, \{2,3\}, \{7\},\{4,8\}\rangle$ of $\{2, \ldots, 8\}$ into 4 parts. Let $\pi_f = 
\langle \{4,5\}, \{1,2\}, \{6\},\{3,7\}\rangle$ be the ordered set partition 
of $[7]$ into 4 parts which results by replacing each element 
$i$ in $\pi^+_f$ by $i-1$.  Then $\pi_f$ is specified 
by the underlying unordered set partition 
$\Gamma_f = \langle \{1,2\},\{3,7\},\{4,5\},\{6\}\rangle$ and 
the permutation $\sigma_f = 3\ 1 \ 4 \ 2$. Finally 
$\sigma_f$ associated to insertion order sequence $h_{\sigma_f} =(2,0,1,0)$ 
and $\Gamma$ is associated to the restricted growth function 
$s_{\Gamma_f} = (0,0,1,2,2,3,1)$.  Thus 
$$Seq(f) = \langle (9,8,6,5,3), (2,0,1,0),(0,0,1,2,2,3,1)\rangle.$$

In the next subsection, we shall provide the basic lemmas about 
ranking and unranking leaves of trees which will allow us to 
reduce the problem of ranking 
$Seq(F_{n,k}) = \{Seq(f):f \in F_{n,k}\}$ according to lexicographic 
order to the problems 
of ranking and unranking decreasing functions according to lexicographic order,
of ranking and unranking insertion sequences of permutations according to lexicographic order, and of ranking and unranking restricted growth functions 
according to lexicographic order. 

\subsection{Basic Lemmas for Ranking and Unranking Algorithms.}

To develop our ranking and unranking algorithms for $\vec{C}_{n,1}^k$ and 
$\vec{C}_{n,1}^{k,L}$, we first need to make some general remarks about ranking and unranking paths in planar trees. Given a rooted planar tree $T$, let $L(T)$ be the numbers of leaves of 
$T$ and $Path(T)$ be the set of paths which go from the root to a leaf.  
Then for any path $p \in Path(T)$, we define the rank of $p$ relative to $T$, 
$rank_T(p)$, to be the number of leaves of $T$ that lie to the left of $p$.

Given two rooted planar trees $T_1$ and $T_2$, $T_1 \otimes T_2$ is the 
tree that results from $T_1$ by replacing each leaf of $T_1$ 
by a copy of $T_2$, see Figure 3. If the vertices of $T_1$ and $T_2$ 
are labeled, then we shall label the vertices of $T_1 \otimes T_2$ according 
to the convention that each vertex $v$ in $T_1$ have the same label in 
$T_1 \otimes T_2$ that it has in $T_1$ and each vertex $w$ in a copy of 
$T_2$ that is decendent from a leaf labeled $l$ in $T_1$ has a label 
$(l,s)$ where $s$ is the label of $w$ in $T_2$. 
Given rooted planar trees $T_1, \ldots, T_k$ where $k \geq 3$, we can define 
$T_1 \otimes T_2 \otimes \cdots \otimes T_k$ by induction as 
$(T_1 \otimes \cdots \otimes T_{k-1}) \otimes T_k$. Similarly if 
$T_1 \ldots, T_k$ are labeled rooted planar trees, we can define 
the labeling of $T_1 \otimes T_2 \otimes \cdots \otimes T_k$ by the same 
inductive process.

\fig{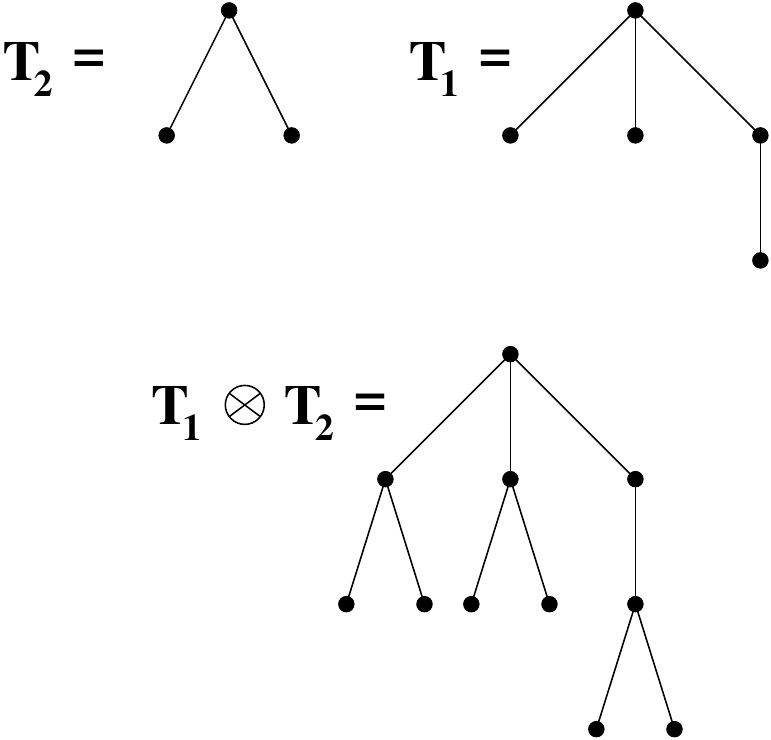}{The operation $T_1 \otimes T_2$.}


Now suppose that we are given two rooted planar trees $T_1$ and $T_2$ and 
suppose that $p_1 \in Path(T_1)$ and $p_2 \in Path(T_2)$.  Then we 
define the path $p_1 \otimes p_2$ in $T_1 \otimes T_2$ which follows 
$p_1$ to its leaf $l$ in $T_1$ and then follows $p_2$ in the copy of 
$T_2$ that sits below leaf $l$ to a leaf $(l,l')$ in $T_1 \otimes T_2$. 
Similarly, given paths $p_i \in T_i$ for $i = 1, \ldots k$, 
we can define a path  $p = p_1 \otimes \cdots \otimes p_k \in Path(T_1 \otimes T_2 \otimes \cdots \otimes T_k)$ by induction as 
$(p_1 \otimes \cdots \otimes p_{k-1}) \otimes p_k$.

Next we give two simple lemmas that tell us how to rank and unrank the set 
of paths in such trees.

\begin{lemma}\label{prodrank}  Suppose that $T_1, \ldots ,T_k$ are rooted planar trees and $T = T_1 \otimes T_2 \otimes \cdots \otimes T_k$.  Then for any path 
$p = p_1 \otimes \cdots \otimes p_k \in Path(T)$, 
\begin{equation}\label{prod-rank}
rank_T(p) = \sum_{j=1}^k rank_{T_j}(p_j) \prod_{l=j+1}^k L(T_{l})
\end{equation} 
\end{lemma}
{\em Proof}.  We proceed by induction on $k$. Let us assume that 
$T_1, \ldots T_k$ are labeled rooted planar trees. 

First suppose that $k =2$ and that $p_1$ is a path that goes from 
the root of $T_1$ to a leaf labeled $1_1$ and $p_2$ goes from the root 
of $T_2$ to a leaf labeled $l_2$.  Thus $p_1 \otimes p_2$ goes from 
the root of $T_1 \otimes T_2$ to the leaf $l_1$ in $T_1$ and then proceeds 
to the leaf $(l_1,l_2)$ in $T_1 \otimes T_2$. Now for each 
leaf $l'$ to the left of $l_1$ in $T_2$, there are $L(T_2)$ leaves 
of $T_1 \otimes T_2$ that lie to left of $(l_1,l_2)$ coming from 
the leaves of the copy of $T_2$ that sits below $l'$. Thus 
there are a total of $L(T_2) \cdot rank_{T_1}(p_1)$ such leaves. The only 
other leaves of $T_1 \otimes T_2$ that lie to left of $p_1 \otimes p_2$ 
are the leaves of the form $(l_1,l'')$ where $l''$ is to left of $p_2$ in 
$T_2$.  There are $rank_{T_2}(p_2)$ such leaves.  Thus there are a total of 
$rank_{T_2}(p_2) + L(T_2) \cdot rank_{T_1}(p_1)$ leaves to left of 
$p_1 \otimes p_2$ and hence 
$$rank_{T_1 \otimes T_2}(p_1 \otimes p_2) = 
rank_{T_2}(p_2) + L(T_2) \cdot rank_{T_1}(p_1) $$
as desired. 

Next assume that (\ref{prod-rank}) holds for $k < n$ and that $n \geq 3$. Then 
\begin{eqnarray*}
&&rank_{T_1 \otimes \cdots \otimes T_n}(p_1 \otimes \cdots \otimes p_n) = \\
&&rank_{(T_1 \otimes \cdots \otimes T_{n-1}) \otimes T_n}
(p_1 \otimes \cdots \otimes p_{n-1}) \otimes p_n) = \\
&& rank_{T_n}(p_n) + L(T_n) 
(\sum_{j=1}^{n-1} rank_{T_j}(p_j) \prod_{l=j+1}^{n-1} L(T_{l})) =\\
&&\sum_{j=1}^n rank_{T_j}(p_j) \prod_{l=j+1}^n L(T_{l}).
\end{eqnarray*}
\hfil$\Box$.

This given, it is easy to develop an algorithm for unranking in a product 
of trees.  The proof of this lemma can be found in \cite{W}. 

\begin{lemma}\label{produnrank}  
Suppose that $T_1, \ldots ,T_k$ are rooted planar trees and $T = T_1 \otimes T_2 \otimes \cdots \otimes T_k$.  Then given a  
$p \in Path(T)$ such that $rank_T(p) = r_0$, 
$p = p_1 \otimes \cdots \otimes p_k \in Path(T)$ where 
$rank_{T_i}(p_i) = q_i$ and 
\begin{eqnarray}
r_0 &=& q_1 \prod_{l=2}^k L(T_{l}) + r_1 \ 
\mbox{where $0 \leq r_1 <  \prod_{l=2}^k L(T_{l})$}, \\
r_1 &=& q_2 \prod_{l=3}^k L(T_{l}) + r_2 \ \mbox{where $0 \leq r_2 <  
\prod_{l=3}^k L(T_{l})$},\\
&\vdots& \\
r_{k-2} &=& q_{k-1}L(T_k) + r_{k-1} \ \mbox{where $0 \leq r_{k-1} < L(T_k)$} \ \mbox{and} \\
r_{k-1} &=& q_k 
\end{eqnarray} 
\end{lemma}

Our idea is to construct trees $T_{{\cal DF}_{n,k}}$,  $T_{{\cal IO}_n}$ and  
$T_{{\cal RG}_{n,k}}$ so that 
\begin{enumerate}
\item the paths of $T_{{\cal DF}_{n,k}}$ correspond to 
the decreasing functions 
in ${\cal DF}_{n,k}$ ranked according to the lexicographic order, 

\item the paths of $T_{{\cal IO}_n}$ correspond to the permutations of $S_n$ 
ranked according to the lexicographic order on their insertion sequences, and 

\item the paths of $T_{{\cal RG}_{n,k}}$ correspond 
to set  restricted growth functions ${\cal RG}_{n,k}$
ranked according to the lexicographic order.
\end{enumerate}

It will then easily follow that the paths of 
$$T_{{\cal DF}_{n,k}} \otimes T_{{\cal IO}_k} \otimes T_{{\cal RG}_{n-2,k}}$$
naturally correspond to the sequences in 
$Seq(F_{n,k}) = \{Seq(f):f \in F_{n,k}\}$ ranked according to lexicographic 
order.  Thus we can use Lemmas \ref{prodrank} and \ref{produnrank} to obtain a ranking and unranking algorithms to 
$Seq(F_{n,k})$ relative to the lexicographic order once we have constructed the trees $T_{{\cal DF}_{n,k}}$,  $T_{{\cal IO}_n}$ and  
$T_{{\cal RG}_{n,k}}$ and developed ranking and unranking algorithms for them. Our next three subsections will be devoted to constructing the trees 
$T_{{\cal DF}_{n,k}}$,  $T_{{\cal IO}_n}$ and  
$T_{{\cal RG}_{n,k}}$ and specifying ranking and unranking algorithms for them. 

\subsection{Ranking and Unranking Decreasing Functions.}

In this subsection, we consider the problem of ranking and unranking the set ${\cal DF}_{n,k}$ of decreasing functions $f:\{1, \ldots, k\} \rightarrow \{1, \ldots, n\}$ relative to lexicographic order. A number of authors have developed ranking and unranking algorithms for ${\cal DF}_{n,k}$.  We shall follow the method of Williamson \cite{W}. First, we identify a function 
$f:\{1, \ldots, k\} \rightarrow \{1, \ldots, n\}$ with the decreasing sequence 
$\langle f(1), \ldots, f(k) \rangle$ where 
$n \geq f(1) > \ldots > f(k) \geq 1$.  We can then think of the sequences as specifying a path in a planar tree $T_{{\cal DF}_{n,k}}$ which can be constructed recursively as follows. At level 1, 
the nodes of $T_{{\cal DF}_{n,k}}$ are labeled   
$k, \ldots, n$ from left to right specifying the choices for $f(1)$. Next below a node $j$ at level one, we attach a tree corresponding to 
$T_{{\cal DF}_{j-1,k-1}}$ where a tree  $T_{{\cal DF}_{a,1}}$ consists of a tree with a single vertex labeled $a$. Figure 3 pictures the tree 
$T_{{\cal DF}_{6,3}}$.

\fig{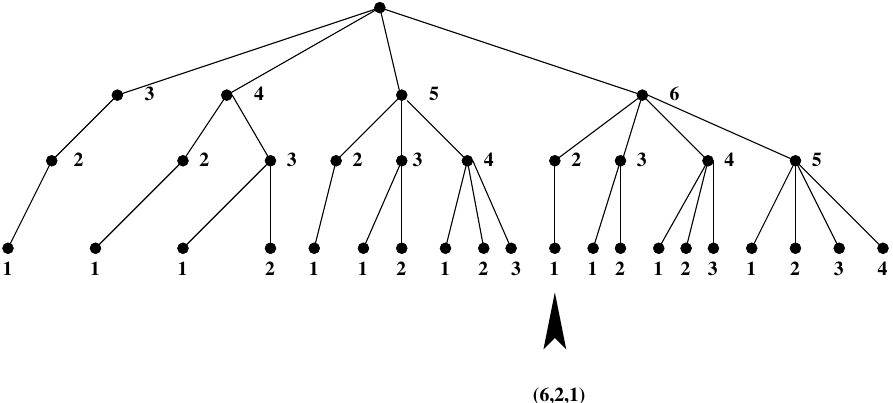}{The tree $T_{{\cal DF}_{7,3}}$.}


Then the decreasing sequence (6,2,1) corresponds to the path from 
the root to the node which is specified with an arrow.  It is clear that the sequences corresponding to the 
paths the tree $T_{{\cal DF}_{6,3}}$ appear in lexicographic order from left to right.  Thus the rank of any sequence $\langle 
f(1), \ldots, f(k) \rangle \in {\cal DF}_{n,k}$ is the number of leaves of 
the tree to the left of the path  corresponding to $\langle 
f(1), \ldots, f(k) \rangle$. Hence the sequence $\langle 6,2,1\rangle$
 has rank 10 in the tree $T_{{\cal DF}_{6,3}}$. 

This given, suppose we are given a sequence $\langle f(1), \ldots, f(k) \rangle$ in $T_{{\cal DF}_{n,k}}$. Then the number of leaves in the subtrees 
corresponding the nodes $k, \ldots, f(1) -1$ are respectively 
$\binom{k-1}{k-1}, \binom{k}{k-1}, \ldots, \binom{f(1)-2}{k-1}$. Thus the total number of leaves in those subtrees is 
$$\binom{k-1}{k-1} + \binom{k}{k-1}+ \cdots + \binom{f(1)-2}{k-1} = 
\binom{f(1)-1}{k}.$$
Here we have used the well known identity that 
$\sum_{s=k-1}^{t-1} \binom{s}{k-1} = \binom{t}{k}$. 
It follows that the rank of 
$\langle f(1), \ldots, f(k) \rangle$ in $T_{{\cal DF}_{n,k}}$ equals 
$\binom{f(1)-1}{k}$ plus the rank of $\langle f(2), \ldots, f(k)\rangle$ in 
$T_{{\cal DF}_{f(1)-1,k-1}}$.  The following result, stated in 
\cite{W}, then easily follows by induction.

\begin{theorem}\label{rankdecfn}
Let $f:\{1, \ldots, k\} \rightarrow \{1, \ldots, n\}$ be a descreasing function. Then the rank of $f$ relative to the lexicographic order on ${\cal DF}_{n,k}$ is 
\begin{equation}\label{decfnrank}
rank_{{\cal DF}_{n,k}}(f) = \binom{f(1)-1}{k} + \binom{f(2)-1}{k-1} + \cdots + \binom{f(k)-1}{1}.
\end{equation}
\end{theorem}

It is then easy to see from Theorem \ref{rankdecfn}, that the following 
procedure, as described by Williamson in \cite{W}, gives the unranking procedure for ${\cal DF}_{n,k}$.\\
\begin{theorem}\label{unrankdecfn}
The following procedure $UNRANK(m)$ computes 
$$f =\langle f(1),\ldots, f(k)\rangle$$
 such that $Rank_{{\cal DF}_{n,k}}(f) =m$ for any 
$1 \leq k \leq n$ and $0 \leq m \leq \binom{n}{k}-1$.\\
\ \\
{\bf Procedure} UNRANK(m) \\
\begin{description}
\item[initialize] $m':=m$, $t:=1$, $s:= k$; {\rm (}$1 \leq k \leq n$, $0 \leq m \leq \binom{n}{k}-1${\rm )}
\item[while] $t \leq k$ \ {\bf do}
\begin{description}
\item[begin]
\begin{description}
\item $f(t)-1 = max\{y: \binom{y}{s} \leq m'\}$;
\item $m' := m' - \binom{f(t)-1}{s}$;
\item $t := t+1$;
\item $s := s-1$;
\end{description}
\item[end]
\end{description}
\end{description}

\end{theorem}

\subsection{Ranking and Unranking Permutations}

The problem of ranking and unranking permutations according 
to lexicographic order on insertion sequences is quite easy. 
We can then think of the insertion sequence $h_\sigma = (h(1), \ldots, 
h(n))$ of a 
permutation $\sigma \in S_n$ as specifying a path 
in a planar tree $T_{{\cal IO}_n}$ which can be constructed recursively as follows. 
At level 1, 
the nodes of $T_{{\cal IO}_n}$ are labeled   
$0, \ldots, n-1$ from left to right specifying the choices for $h_1$. Next below a node $j$ at level one, we attach a tree corresponding to 
$T_{{\cal IO}_{n-1}}$ where a tree  $T_{{\cal IO}_1}$ consists of a tree with a single vertex labeled $0$. Figure \ref{figure:TS3} pictures the tree 
$T_{S_3}$.

\fig{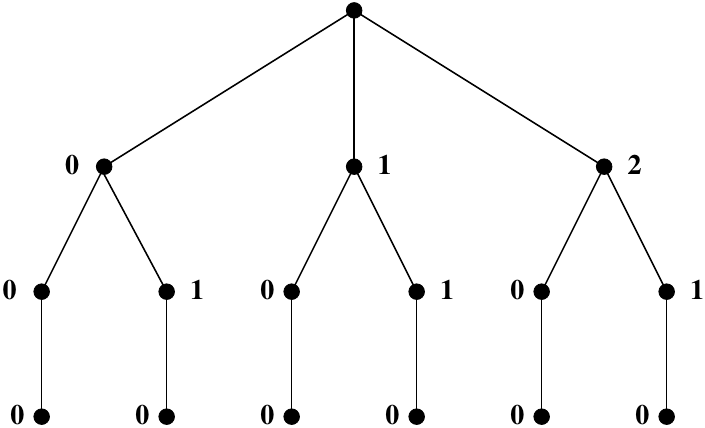}{The tree $T_{{\cal IO}_3}$.}


This given, it easy to see from our definitions that the following theorems,
which  can be found in \cite{W}, hold.

\begin{theorem}\label{thm:permrank}
Let $S_n$ denote the symmetric group of all permutations of $[n]$.  
Order $S_n$ by defining $\sigma < \tau$ iff $h_\sigma \leq_{lex} h_\tau$, 
then the rank of $\sigma$, $Rank(\sigma) = |\{\tau \in S_n: \tau < \sigma\}|$, can be computed by the formula
\begin{equation}\label{eq:permrank}
Rank(\sigma) = \sum_{k=1}^{n-1} h(i) \cdot (n-i)!
\end{equation}
if $h_\sigma = (h(1), \ldots, h(n))$. 
\end{theorem}

\begin{theorem}\label{thm:permunrank}
Let $S_n$ denote the symmetric group of all permutations of $[n]$.  
Order $S_n$ by defining $\sigma < \tau$ iff $h_\sigma \leq_{lex} h_\tau$. 
Then if $0 \leq p \leq n!-1$, then we construct the insertion sequence 
$h_\sigma = (h(1), \ldots, h(n))$ of the permuation $\sigma \in S_n$ such 
that $Rank(\sigma) = p$ as follows. Let $p = p_0$, then set $h_n =0$ and 
compute the following $k = 0,  \ldots, n-2$:
\begin{eqnarray*}
p_0 &=& h(1) (n-1)! + p_1 \ \mbox{where $p_1 < (n-1)!$}\\
p_1 &=& h(2) (n-2)! + p_2 \ \mbox{where $p_2 < (n-2)!$}\\
&\vdots&\\
p_k &=& h(k+1) (n-k)! + p_k \ \mbox{where $p_k < (n-k)!$}\\  
&\vdots&\\
\end{eqnarray*} 
\end{theorem}

\section{Ranking and Unranking Algorithms for Unordered Set Partitions.}

The main purpose of this section is to give ranking and unranking algorithms
for the set of all unordered set partitions of $[n]$ into $k$ parts,
${\cal S}_{n,k}$. We shall follow the general approach of Williamson
\cite{W}
and give algorithms to rank and unrank the set of surjective restricted growth
functions relative to lexicographic order.
A restricted growth function
$f:[n] \rightarrow \{0, 1, \ldots, k-1\}$ is a  function such that
\begin{description}
\item[(i)] $f(1)= 0$ and
\item[(ii)] for all $1 < i \leq n$,
$f(i) \leq  1 + max(\{f(1), \ldots, f(i-1)\})$.
\end{description}
We let ${\cal RG}_{n,k}$ denote the set of all restricted growth functions $f$
that map  $[n]$ {\bf onto\/} $\{0, 1, \ldots, k-1\}$.
There is a natural one-to-one correspondence $I$ between ${\cal RG}_{n,k}$ and
${\cal S}_{n,k}$.
That is, suppose $\pi = \langle \pi_1, \ldots, \pi_k \rangle$ is a set
partition
of $[n]$ into $k$ parts where $min(\pi_1) < \ldots < min(\pi_k)$.
Then define
$I(\pi) = f_{\pi} \in {\cal RG}_{n,k}$ by letting
$f_\pi(i) = j-1$ where $i \in \pi_j$ for all $i \in [n]$.
It is easy to see that  if $f \in {\cal RG}_{n,k}$, then
$I^{-1}(f) = \pi_f = \langle f^{-1}(0), \ldots, f^{-1}(k-1)\rangle$.
We shall identify a restricted growth function $f \in {\cal RG}_{n,k}$ with the
sequence $\langle f(1), \ldots, f(n)\rangle$ and then order ${\cal RG}_{n,k}$
via the lexicographic order on such sequences.

Let $\langle f(1), \ldots, f(n-j)\rangle$, $0\leq j\leq n-1$, be a
restricted growth function and suppose that 
$max(\{f(1), \ldots, f(n-j)\}) = m$ where  $0\leq m\leq k-1$.
Let
\begin{eqnarray}
&&{\cal E}^{(k)}(j,m,f(1), \ldots, f(n-j))= \\
&&\{
\langle x_{n-j+1}, \ldots, x_n\rangle\mid
\langle f(1), \ldots, f(n-j), x_{n-j+1}, \ldots, x_n \rangle\in {\cal RG}_{n,k}
\}. \nonumber
\end{eqnarray}
The set ${\cal E}^{(k)}(j,m,f(1), \ldots, f(n-j))$ represents all ways of
extending the sequence \\
$f(1), \ldots, f(n-j)$ to a surjective  restricted growth function. 
Note that \\
(1) ${\cal E}^{(k)}(0,m,f(1), \ldots, f(n)) = \emptyset$ for $m\neq k-1$,\\
(2) ${\cal E}^{(k)}(0,k-1, f(1), \ldots, f(n)) = \{\epsilon\}$ where $\epsilon$ is the empty string, and \\
(3) for $0<j\leq n-1$,
\begin{eqnarray*}
&&{\cal E}^{(k)}(j,m,f(1), \ldots, f(n-j)) = \\
&&\bigcup_{t=0}^{m} t*{\cal E}^{(k)}(j-1,m,f(1), \ldots, f(n-j), t) \\
&&\cup 
(m+1)*{\cal E}^{(k)}(j-1,m+1,f(1), \ldots, f(n-j), m+1).
\end{eqnarray*}
This union is disjoint.
The notation $t*X$ means concatenate $t$ with each element in $X$.
Recall that we have assumed above that
$max(\{f(1), \ldots, f(n-j)\}) = m$.
In the sets
$${\cal E}^{(k)}(j-1,m,f(1), \ldots, f(n-j), t)$$
we have $0\leq t \leq m$,
thus in the function
$$\langle f(1), \ldots, f(n-j), f(n-j+1)\rangle =
\langle f(1), \ldots, f(n-j), t\rangle
$$
we still have $max(\{f(1), \ldots, f(n-j+1)\}) = m$.
Given that $max(\{f(1), \ldots, f(n-j)\}) = m$, we
can also assign $f(n-j+1) = m+1$ and still have
$\langle f(1), \ldots, f(n-j+1)\rangle$
a restricted growth function.
The set
$(m+1)*{\cal E}^{(k)}(j-1,m+1,f(1), \ldots, f(n-j), m+1)$
includes all such elements in
${\cal E}^{(k)}(j,m, f(1), \ldots, f(n-j))$.
Note that ${\cal E}^{(k)}(j,m, f(1), \ldots, f(n-j))$ does
not depend on the actual values of $f(1), \ldots, f(n-j)$.
Setting
$$|{\cal E}^{(k)}(j,m, f(1), \ldots, f(n-j))| = E^{(k)}(j,m)$$
we have the following recursion
$$E^{(k)}(j,m) = (m+1)E^{(k)}(j-1,m) + E^{(k)}(j-1,m+1)$$
The initial values are $E(0,m)=0$ if $m\neq k-1$, $E(0,k-1)=1$.
This recursion will be the basis for our ranking and unranking
procedures for ${\cal RG}_{n,k}$ and hence for ${\cal S}_{n,k}$.

For example, Table 3 below lists the values of $E^{(4)}(m,n)$. 

\begin{center}
{\normalsize 
\begin{tabular}{|l||l|l|l|l|l|l|l|l|}
\hline 
 $\begin{matrix} \ & m \\ n & \ \end{matrix}$& 0 \ \ \ & 1 \ \ \ & 2 \ \ \  &3 \ \ \  & 4 \ \ \ & 5 \ \ \ & 6 \ \ \ & 7\ \ \ \\
\hline \hline
0&0&0&0&1&0&0&0&0\\
\hline 
1&0&0&1&4&0&0&0& \\
\hline
2&0&1&7&16&0&0&&\\
\hline
3&1&9&37&64&0&&&\\
\hline
4&10&55&175&256&&&&\\
\hline
5&65&285&781&&&&&\\
\hline 
6&350&1351&&&&&&\\
\hline
7&1701&&&&&&&\\
\hline
\end{tabular}}
\end{center}
\begin{center}
{\bf Table 3} $E^{(4)}(m,n)$
\end{center}

It is easy to construct a planar 
tree $T_n$ whose paths correspond to restricted 
growth functions $(\alpha_{n-1}, \alpha_{n-2}, \ldots, \alpha_0)$ 
by having the root of the $T_n$ labelled with $0$ and by 
induction if a node corresponds to a path 
$(\alpha_{n-1}, \ldots, \alpha_{t+1})$, then the 
descendants of that node would be labeled from left to right with $0, 
\ldots, m_t, m_t+1$ where $m_t = max(\{\alpha_{n-1}, \ldots, \alpha_{t+1}\})$; 
see Figure ref{figure:RecT}. The tree $T_{{\cal RG}_{n,k}}$ would consists of those nodes which lies on a path 
$(\alpha_{n-1}, \alpha_{n-2}, \ldots, \alpha_0)$ such that 
$max(\{\alpha_{n-1}, \alpha_{n-2}, \ldots, \alpha_0\}) = k-1$.  \
For example, $S_{4,3} =6$ and the tree 
$T_{{\cal RG}_{4,3}}$ is pictured in Figure \ref{figure:TS43}

\fig{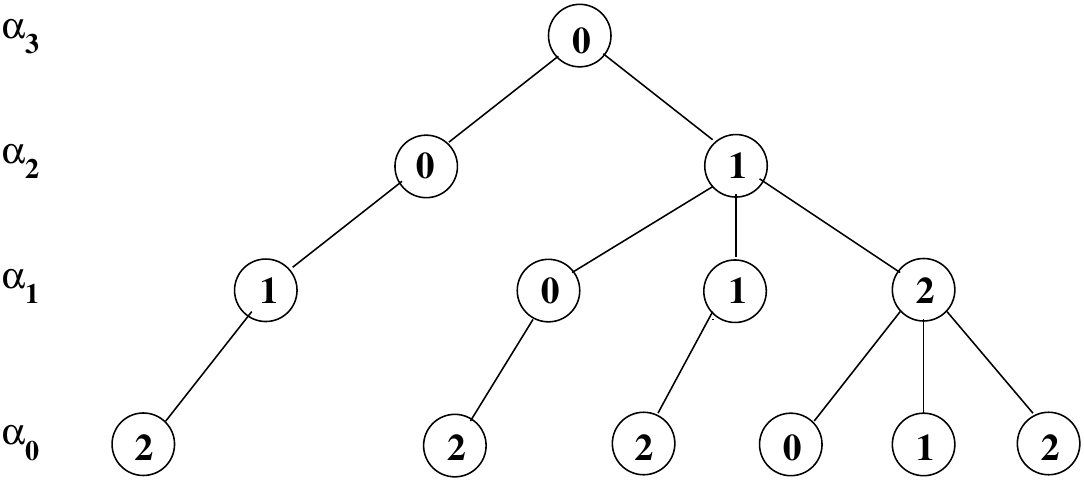}{The tree $T_{{\cal RG}_{4,3}}$.}

More generally, it is easy to see that if we are at node corresponding 
to the sequence $(\alpha_{n-1}, \ldots, \alpha_{t+1})$ in 
$T_{{\cal RG}_{n,k}}$ where 
$m_t = max(\alpha_{n-1}, \ldots, \alpha_{t+1}) < k-1$, then, as pictured 
in figure \ref{figure:RecT}, our choices for $\alpha_t$ is either $0,1, \ldots, m_t, m_t+1$. Moreover by our definition of $E^{(k)}(t,m_t)$, the number of 
leaves in each of the subtrees where we branch to either $0,1, \ldots, m_t$ is exactly $E^{(k)}(t,m_t)$ and the number of leaves in the subtree where we 
branch to $m_t+1$ is $E^{(k)}(t,m_t+1)$. It follows that if we take the branch corresponding to $j$ out the node corresponding to the sequence $(\alpha_{n-1}, \ldots, \alpha_{t+1})$, then the number of nodes in the subtrees to the left of that branch is $j\cdot E^{(k)}(t,m_t)$.  
If we are at node corresponding 
to the sequence $(\alpha_{n-1}, \ldots, \alpha_{t+1})$ in 
$T_{{\cal RG}_{n,k}}$ where 
$m_t = max(\alpha_{n-1}, \ldots, \alpha_{t+1}) = k-1$, then we can only 
branch to either $0, 1, \ldots, m_t$ but it is still the case 
that, if we take the branch corresponding to $j$ out the node corresponding to the sequence $(\alpha_{n-1}, \ldots, \alpha_{t+1})$, then the number of nodes in the subtrees to the left of that branch is $j\cdot E^{(k)}(t,m_t)$.   
The following result of Williamson \cite{W} for ranking sequences $f = (\alpha_{n-1}, \alpha_{n-2}, \ldots, \alpha_0)$ be an element of ${\cal RG}_{n,k}$ relative to lexicographic order, then easily follows by induction. 

\fig{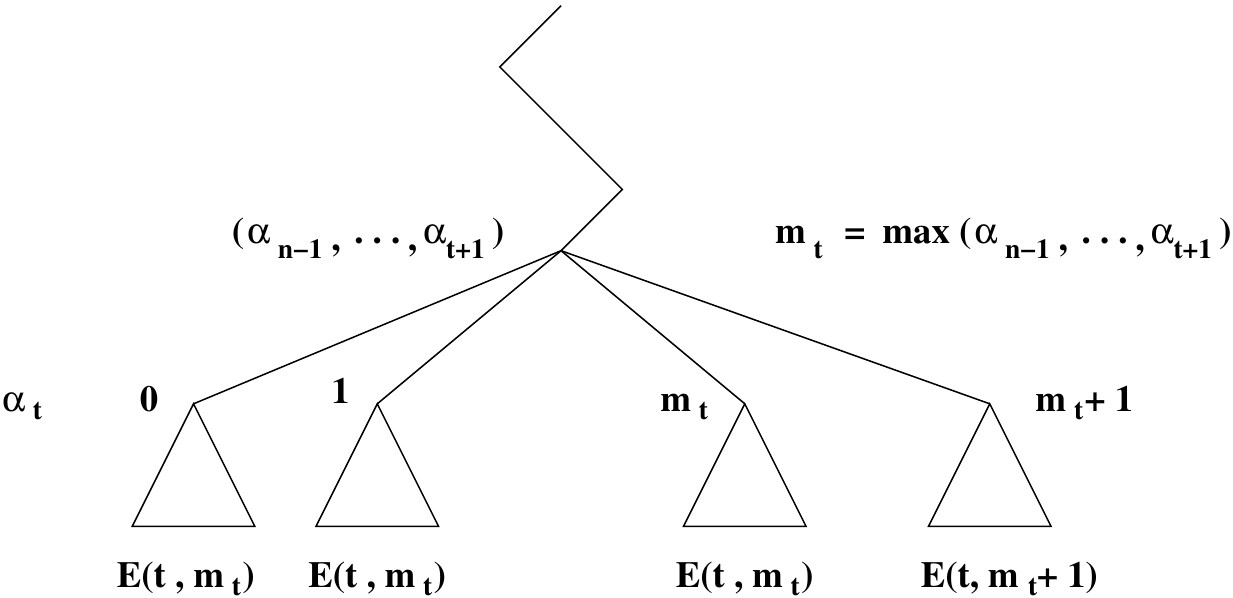}{Typical node in $T_{{\cal RG}_{n,k}}$.}

\begin{theorem} \label{thm:RGrank} 
Let $f = (\alpha_{n-1}, \alpha_{n-2}, \ldots, \alpha_0)$ be an element 
of ${\cal RG}_{n,k}$ and for each $t = 0, 1, \ldots, n-2$, let 
$m_t = max(\{\alpha_{n-1}, \ldots, \alpha_{t+1}\})$.  Then 
\begin{equation}\label{eq:RGrank}
rank_{{\cal RG}_{n,k}}(f) = \sum_{t=0}^{n-2} \alpha_t E^{(k)}(t,m_t).
\end{equation}
\end{theorem}

For example, suppose that 
$f =(0,1,0,2,0,3,1,2) = (\alpha_{7},\ldots, \alpha_0)\in {\cal RG}_{8,4}$.  
Then for $t =0, \ldots, 6$, we let $m_t = 
max(\{\alpha_7, \ldots, \alpha_{t+1}\})$ 
so that $m_0=3, m_1=3, m_2 =2, m_3 =2, m_4 =1, m_5 = 1$ and $m_6 =0$. Thus 
using the values of $E^{(4)}(n,m)$ from Table 3, we see that  
\begin{eqnarray*}
rank_{{\cal RG}_{n,k}}(f) &=& \sum_{t=0}^{6} \alpha_t E^{(4)}(t,m_t) \\
&=& 2 \cdot E^{(4)}(0,3) + 1 \cdot E^{(4)}(1,3) + 3 \cdot E^{(4)}(2,2) + 0 \cdot E^{(4)}(3,2)+ \\
&& 2 \cdot E^{(4)}(4,1) + 0 \cdot E^{(4)}(5,1) + 1 \cdot E^{(4)}(6,0)\\
&=& 2\cdot 1 + 1 \cdot 4 + 3 \cdot 7 + 0 \cdot 37 + \\
&& 2 \cdot 55 + 0 \cdot 285 + 1 \cdot 350 \\
&=& 487.
\end{eqnarray*}

One can also use Figure \ref{figure:RecT} to give an inductive proof of the validity of the following unranking algorithm for the element of ${\cal RG}_{n,k}$ relative to lexicographic order. 

\begin{theorem} \label{thm:RGunrank} Let $0 \leq r < S_{n,k}$ where $S_{n,k}$ is the number of set partitions of $[n]$ into $k$ parts.  Then 
we can construct 
$f = (\alpha_{n-1}, \alpha_{n-2}, \ldots, \alpha_0) \in {\cal RG}_{n,k}$ 
such that $rank_{{\cal RG}_{n,k}}(f) = r$ as follows.
\begin{description}
\item[Step $n-1$] Set $\alpha_{n-1} = 0$.

\item[Step $n-2$] Set $m_{n-2} = 0$ and let $\alpha_{n-2}$ equal the 
largest $s$ such that $0 \leq s \leq m_{n-1} +1$ and 
$s E^{(k)}(n-2,m_{n-1}) \leq r$.  Then set 
$$r : r - \alpha_{n-2}E^{(k)}(n-2,m_{n-1}).$$

\item[Step $t$] Assume that we have defined 
$\alpha_{n-1}, \ldots, \alpha_{t+1}$. Let 
$m_t = max(\{\alpha_{n-1}, \ldots, \alpha_{t+1}\})$.  Then set 
$\alpha_{t}$ equal the 
largest $s$ such that $0 \leq s \leq m_{t} +1$ and 
$s E^{(k)}(t,m_{t}) \leq r$.  Then set 
$$r : r - \alpha_t E^{(k)}(t,m_{t}).$$  
\end{description}
\end{theorem}

Thus, for example, to compute the 
$f = (\alpha_7, \ldots, \alpha_0) \in {\cal RG}_{8,4}$ whose rank is 1000, we would carry out the following steps. \\
\begin{description}
\item[Step 7.] We set $r= 1000$ and $\alpha_7 =0$. 

\item[Step 6.] $m_6 = \alpha_7 = 0$. $E^{(4)}(6,0) = 350$.  Since 
$1\cdot E^{(4)}(6,0) = 350 < 1000$, then $\alpha_6 =1$ and we set 
$r = 1000 -350 =650$. 

\item[Step 5.] $m_5 = max(\{0,1\}) =1$. $E^{(4)}(5,1) = 285$.  Since 
$2\cdot E^{(4)}(5,1) = 570  < 650$, 
then $\alpha_5 =2$ and we set 
$r = 650 -570 =80$. 

\item[Step 4.] $m_4 = max(\{0,1,2\}) =2$. $E^{(4)}(4,2) = 175$.  Since 
$0\cdot E^{(4)}(4,2) = 0  < 80 < 1 \cdot E^{(4)}(4,2) = 175$, 
then $\alpha_4 =0$ and we set 
$r = 80 - 0 \cdot 175  =80$.

\item[Step 3.] $m_3 = max(\{0,1,2,0\}) =2$. $E^{(4)}(3,2) = 37$.  Since 
$2\cdot E^{(4)}(3,2) = 74  < 80 < 3 \cdot E^{(4)}(3,2) = 111$, 
then $\alpha_4 =2$ and we set 
$r = 80 - 2 \cdot 37  =6$.
  
\item[Step 2.] $m_2 = max(\{0,1,2,0,2\}) =2$. $E^{(4)}(2,2) = 7$.  Since 
$0\cdot E^{(4)}(2,2) = 0  < 6 < 1 \cdot E^{(4)}(3,2) = 7$, 
then $\alpha_2 =0$ and we set 
$r = 6 - 0 \cdot 7  =6$.

\item[Step 1.] $m_1 = max(\{0,1,2,0,2,0\}) =2$. $E^{(4)}(1,2) = 1$.  Since 
$3\cdot E^{(4)}(1,2) = 3  < 6$, 
then $\alpha_1 =3$ and we set 
$r = 6 - 3 \cdot 1  =3$.

\item[Step 0.] $m_0 = max(\{0,1,2,0,2,0,3\}) =3$. $E^{(4)}(0,3) = 1$.  Since 
$3\cdot E^{(4)}(0,3) = 3  \leq 3 < 4 \cdot E^{(4)}(0,3) =4$, 
then $\alpha_0 =3$ and we set 
$r = 3 - 3 \cdot 1  =0$.
\end{description}
Thus the restricted growth function $f = (0,1,2,0,2,0,3,3)$ is the element of 
rank 1000 in ${\cal RG}_{8,4}$.

\section{Ranking and Unranking Algorithms for $\vec{C}_{n,1}^{k,L}$ or $\vec{C}_{n,1}^k$}

The main goal of this section is to give our final ranking and unranking algorithms for the sets $\vec{C}_{n,1}^{k,L}$ or $\vec{C}_{n,1}^k$.  

We start by considering 
the ranking and unranking algorithms for $\vec{C}_{n,1}^k$. Recall 
$$|\vec{C}_{n,1}^k| = \binom{n}{k} \times (n-k)! \times S_{n-2,n-k}.$$
\ \\
{\bf Ranking Algorithm for $\vec{C}_{n,1}^k$.}\\
\ \\
{\bf Step 1}  Given $T \in \vec{C}_{n,1}^k$, construct $f = \Theta^{-1} (T)$.\\
We assume that we start with the data structure for 
a tree $T \in \vec{C}_{n,1}^k$ which consists of pairs $\langle i, j\rangle$ such that $(i,j)$ is a directed edge in $T$. We then compute $f= \Theta^{-1}(T)$. 
By our results in section 3, we can construct the set of pairs 
$\langle i,j\rangle$ such that $f(i) =j$ in linear time. \\
\ \\
\ \\
{\bf Step 2.}  Find $g \in {\cal DF}_{n,k}$, $h \in {\cal IO}_{n-k}$, and 
$s \in {\cal RG}_{n-2,n-k}$ such that 
$Seq(f) = \langle (g(1), \ldots g(k)), (h(1), \ldots, h(n-k)),(s(1), \ldots, s(n-2))\rangle$.  \\
Note that by making one pass through the data structure for $f$, we can construct of table $Table(f) = \langle (1,V_1), \ldots, (n,V_n)\rangle$ where for each $i$, $V_i = \emptyset$ if 
$f^{-1}(i) = \emptyset$ or $V_i$ a linked list of the elements of 
$f^{-1}(i)$ in increasing order 
if $f^{-1}(i) \neq \emptyset$. For example if we start 
with the tree $T \in \vec{C}_{9,1}^5$ and $f=\Theta^{-1}(T)$ pictured in Figure 3, then one can see from Table 1 that we would produce the following table.
\begin{equation}\label{ex1}
Table(f) = \langle (1,(5,6)),(2,(2,3)),(3, \emptyset), (4,(7)), (5, \emptyset),(6, \emptyset),(7,(4,8)), (8, \emptyset),(9, \emptyset)\rangle.
\end{equation}
We can then read Table(f) from right to left to construct 
the decreasing sequence $$(g(1), \ldots, g(k))$$ of those $i$ such that $f^{-1}(i) = \emptyset$.  Such elements correspond to the leaves of the tree $T$ by Theorem 2. 
Thus is clear that we can construct $(g(1), \ldots, g(k))$ from $T$ in linear time. For the tree $T$ pictured in Figure 3, we would produce 
\begin{equation}\label{ex2}
(g(1), \ldots, g(5)) = (9,8,6,5,3).
\end{equation}

By reading $Table(f)$ from left to right, we can construct the ordered 
set partition of $\{2, \ldots, n-1\}$ into $n-k$ parts, 
$\pi^+(f) = \langle \pi^+_1, \ldots, \pi^+_{n-k}\rangle$, which consists of  
the $V_i$'s which are nonempty in 
$Table(f)$ and a set of pointer from $\pi^+_j$ to $j$.  For our example, 
\begin{equation}\label{ex3}
\pi^+(f) = \langle (5,6),(2,3),(7),(4,8) \rangle.
\end{equation}
Next we make a pass through $\pi^+(f)$ and create 
\begin{description}
\item[(i)] a new ordered set partition 
$\pi(f) =\langle \pi_1, \ldots, \pi_{n-k}\rangle$ of $[n-2]$ into $n-k$ parts along with pointers from $\pi_j$ to $j$ by  replacing 
each element $i$ by $i-1$,  
\item[(ii)] a sequence 
$u(f) = (u(1), \ldots, u(n-k)$ where 
$u(i) = min(\pi^+_i) -1$, and 
\item[(iii)] a sequence  $v(f) = (v(1), \ldots, v(n-2))$  where 
$v(i) = 0$ if $i$ appears in  $u(f)$ and $v(i) = 1$ otherwise. 
\end{description}
For our given example of $f$, 
we would produce 
\begin{equation}\label{ex4}
\pi(f) = \langle (4,5),(1,2),(6),(3,7) \rangle,
\end{equation}
\begin{equation}\label{ex5}
u(f) =(4,1,6,3),
\end{equation}
and
\begin{equation}\label{ex6}
v(f) = (0,1,0,0,1,0,1).
\end{equation}
Again, it is easy to see that we can construct $\pi(f)$, $u(f)$, and $v(f)$ in linear time from $T$. 

Next we use $v(f)$ to construct $\delta(f) = (\delta(1), \ldots, \delta(n-k))$ where $\delta(1) =0$ and for $1 < i \leq n-2$, $\delta(i) = \sum_{j=1}^{i-1} v(j)$.  Thus $\delta(i)$ is the number of $j < i$ such $j$ is not a 
minimum element in one of 
the parts $\pi(f)$. It follows that if we let $\sigma = (\sigma(1), \ldots, 
\sigma(n-k)) = (u(1) - \delta(u(i)), \ldots, u(n-k) -\delta(n-k))$, then 
$\sigma$ is a permutation of $n-k$ which represent the relative order of parts of $\pi(f)$ according to increasing minimal elements. For example, for our given example of $f$, 
\begin{equation}\label{ex7}
\delta(f) = (0,0,1,1,1,2,2)
 \end{equation}
and 
\begin{equation}\label{ex8}
\sigma = (4 -1,1 -0,6-2,3 -1) = (3,1,4,2)
\end{equation}
Note that if we order the parts of $\pi(f)$ according to increasing 
minimal elements to get 
$\pi^-(f) =\langle \pi^-_1, \ldots, \pi^-_{n-k}\rangle$, the 
$\pi_i = \pi^-_{\sigma(i)}$. For our given example of $f$, 
\begin{equation}\label{ex9}
\pi^-(f) = \langle (1,2),(4,5),(6),(3,7)\rangle.
\end{equation}
Note that we change the pointer on $\pi(f) =\langle \pi_1, \ldots, \pi_{n-k}\rangle$ so that $\pi_j$ pointers to $\sigma_j-1$, the part of $\pi(f)$ which has 
the smallest minimal element will point to $0$, the part of $\pi(f)$ which has 
the second smallest minimal element will point to $1$, etc.  It follows 
that if we make another 
 pass through $\pi(f)$ and as we encounter an element $x$, we 
record $\sigma_j-1$ where the part that contains $x$ in $\pi(f)$ points to $j$, we can construct the  sequence $s(f) = (s(1), \ldots, s(n-k))$, then $s(f)$ will be the restricted growth function associated to the unorderd set partition 
$\pi^-(f)$.  For our given $f$, this process would produce 
\begin{equation}\label{ex10}
s(f) =(0,0,1,2,2,3,1).
\end{equation}
Again it is easy to see that we can produce $s(f)$ in linear time from $T$. 
 
Finally we need to construct the direct insertion sequence 
$h(f) = (h(1), \ldots, h(n-k))$ for $\sigma$. We can produce the insertion sequence by recursion. That is, $h(1) = \sigma(1) -1$ and $(h(2), \ldots, h(n-k))$ is the insertion sequence for the permutation $\tau = (\tau(1), \ldots, \tau(n-k-1))$ where 
$\tau(i) = \sigma(i+1)$ if $\sigma(i+1) < \sigma(1)$ and $\tau(i) = \sigma(i+1) -1$ if $\sigma(i+1) > \sigma(1)$. Thus it requires $0((n-k)^2)$ steps to 
produce $h(f)$.  For our given example of $f$, 
\begin{equation}\label{ex11}
h(f) = (2,0,1,0).
\end{equation}
It follows that it require $O(n^2)$ step to produce $Seq(f)$. \\
\ \\
\ \\
{\bf Step 3}.  Use the algorithms of section 4 to compute 
$rank_{{\cal DF}_{n,k}}(g(f))$, $rank_{{\cal IO}_{n-k}}(h(f))$, and 
$rank_{{\cal RG}_{n-2,n-k}}(s(f))$.\\
It is easy to see from our ranking algorithms of section 4 that 
if we start with table of binomial coefficients and $E^{(n-k)}(r,s)$, then 
it requires $O(n)$ operations of addition, multiplication, and comparisions to compute each of $rank_{{\cal DF}_{n,k}}(g(f))$, 
$rank_{{\cal IO}_{n-k}}(h(f))$, and $rank_{{\cal RG}_{n-2,n-k}}(s(f))$. Since 
the integers involved have length at most $O(nlog(n))$, it follows that it requires $O(n^2log(n))$ steps to compute $rank_{{\cal DF}_{n,k}}(g(f))$, $rank_{{\cal IO}_{n-k}}(h(f))$, and 
$rank_{{\cal RG}_{n-2,n-k}}(s(f))$.

For our given $f$, 
\begin{eqnarray}\label{ex12}
rank_{{\cal DF}_{9,5}}(g(f)) &=& \sum_{i=1}^5 \binom{g(i) -i}{6-i} \\
&=& \binom95 + \binom74 + \binom53 + \binom42 + \binom21 \nonumber \\
&=& 56 + 35+ 10 + 6 +2 = 109, \nonumber
\end{eqnarray}
\begin{eqnarray}\label{ex13}
rank_{{\cal IO}_{4}}(h(f)) &=& \sum_{i=1}^4 h(i) (n-i)! \\
&=& (2 \times 3!) + (0 \times 2!) + (1 \times 1!) + (0 \times 0!)  \nonumber \\
&=& 7, \nonumber
\end{eqnarray}
and if $s(f) = (\alpha_6, \alpha_5, \ldots, \alpha_0)$ and 
$m_i = max(\{\alpha_6, \ldots, \alpha_{i+1}\})$, then  
\begin{eqnarray}\label{ex14}
&&rank_{{\cal RG}_{7,4}}(s(f)) =  \sum_{i=0}^5  \alpha_i E^{(4)}(i,m_i) \\
&&= (0 \times E^{(4)}(5,0)) + (1 \times E^{(4)}(4,0))+ 
(2 \times E^{(4)}(3,1)) \nonumber \\
&&+ ( 2 \times E^{(4)}(2,2)) + 
(3 \times E^{(4)}(1,2)) + (1 \times E^{(4)}(0,3))\nonumber \\
&&= (0 \times 65) + (1 \times 10) + (2 \times 9) + (2 \times 7) + (3 \times 1) + (1 \times 1) = 46. \nonumber
\end{eqnarray}
\ \\
\ \\
{Step 4} It then follows from our analysis in section 4 that 
\begin{eqnarray}\label{ex15}
rank_{\vec{C}_{n,1}^k} &=& 
rank_{{\cal DF}_{n,k}}(g(f)) \cdot (n-k)! S_{n-2,n-k} \\
&&+ rank_{{\cal IO}_{n-k}}(h(f)) \cdot S_{n-2,n-k} \nonumber \\
&& + rank_{{\cal RG}_{n-2,n-k}}(s(f)). \nonumber
\end{eqnarray}
For our given $f$, $n=9$, $k =5$, $(n-k)! = 4! =24$, $S_{7,4} = 350$, $4! \times 350 = 8400$. Thus 
\begin{equation}\label{ex16}
rank_{\vec{C}_{9,1}^5}(T) = (109 \times 8400) + (7 \times 350) + (46) = 918,096
\end{equation}
Thus the $T$ pictured in Figure 3 is the tree with rank 918,096 in 
$\vec{C}_{9,1}^5$.

The unranking algorithm for $\vec{C}_{n,1}^k$ is relatively straight forward given the unranking algorithms discussed in section 4.\\
\ \\
\ \\
{\bf Unranking Algorithm for $\vec{C}_{n,1}^k$.}\\
Problem:  Given $r$ such that $0 \leq r < \binom{n}{k}\times (n-k)! \times S_{n-2,n-k}$, find $T \in \vec{C}_{n,1}^k$ such that 
$rank_{\vec{C}_{n,1}^k}(T) =r$. \\
\ \\ 
{\bf Step I} Find $r_1$ $r_2$, and $r_3$ such that 
\begin{eqnarray}\label{EX1}
r &=& r_1 \times  (n-k)! S_{n-2,n-k} + u_1 \ \mbox{where} \ 0 \leq u_1  < (n-k)! S_{n-2,n-k} \nonumber \\
u_1 &=& r_2 \times S_{n-2,n-k} + u_2 \ \mbox{where} \ 0 \leq u_2 < S_{n-2,n-k} \nonumber \\
r_3 &=& u_2.
\end{eqnarray}
\ \\
{\bf Step II} Find $g \in {\cal DF}_{n,k}$, $h \in {\cal IO}_{n-k}$ and 
$s \in {\cal RG}_{n-2,n-k}$ such that 
\begin{eqnarray}\label{EX2}
r_1 &=& rank_{{\cal DF}_{n,k}}(g) \\
r_2 &=& rank_{{\cal IO}_{n-k}}(h) \\
r_3 &=& rank_{{\cal RG}_{n-2,n-k}}(s) 
\end{eqnarray}
\ \\
{\bf Step III} Find $f \in {\cal F}_n$ such that 
\begin{equation}\label{EX3}
Seq(f) = \langle (g(1), \ldots, g(k)), (h(1), \ldots, h(n-k)), (s(1), \ldots, s(n-2))\rangle.
\end{equation}
\ \\
{\bf Step IV} Let $T = \Theta(f)$.\\
\ \\
Note that since we are dealing with numbers whose length is of $O(nlog(n)$, Step I requires $O(nlog(n))$ steps.  The unranking algorithms for ${\cal DF}_{n,k}$, ${\cal IO}_{n-k}$ and ${\cal RG}_{n-2,n-k}$ each require 
$O(n)$ operations of addition, mulitiplication, division and comparisons. 
Again, since we are dealing with numbers whose length is of $O(nlog(n)$, Step II requires $O(n^2log(n))$ steps.  It is not difficult to see that we can reconstruction $f$ from $g$, $h$ and $s$ in $O(n^2)$ so that Step III requires $O(n^2)$ steps. Finally Step IV can be carried out in $O(n)$ steps so that the unranking proceedure requires $O(n^2log(n))$ steps. 

For an example of our unranking algorithms for $\vec{C}_{n,1}^k$, suppose 
that we want to find the tree $T \in \vec{C}_{9,1}^5$ whose rank is 600,000. 
In this case, $(9-5)!S_{7,4} = 8400$ and $S_{7,4} = 350$.  Thus for Step I, 
we find that 
\begin{eqnarray*}\label{EX4}
600,000 &=& 71 \times 8400 + 3600 \\\
3600 &=& 10 \times 350 + 100.
\end{eqnarray*}
Thus $r_1 = 71$, $r_2 = 10$ and $r_3 =100$.

For Step II, first we apply the uranking proceedure for ${\cal DF}_{9,5}$ 
from section 4.2 to 
find $g \in {\cal DF}_{9,5}$ such that $rank_{{\cal DF}_{9,5}}(g) = 71$. 
We set $m' =71$.\\
\ \\
Then $g(1) -1 = max\{y: \binom{y}{5} \leq 71\}$.  Note that 
$ \binom85 =56 < 71 < 126 = \binom95$. Thus $g(1) -1 = 8$ and hence 
$g(1) =9$. We then set $m' = 71-56 =15$.\\
\ \\
Then $g(2) -1 = max\{y: \binom{y}{4} \leq 15\}$.  Note that 
$\binom64 =15$. Thus $g(2) -1 = 6$ and hence 
$g(2) =7$. We then set $m' = 15 -15 = 0$.\\
\ \\
Then $g(3) -1 = max\{y: \binom{y}{3} \leq 0\}$.  Note that 
$\binom23 =0$. Thus $g(3) -1 = 2$ and hence 
$g(3) =3$. We then set $m' = 0$.\\
\ \\
Then $g(2) -1 = max\{y: \binom{y}{2} \leq 0\}$.  Note that 
$\binom12 =0$. Thus $g(2) -1 = 1$ and hence 
$g(2) =2$. We then set $m' = 0$.\\  
\ \\
Then $g(1) -1 = max\{y: \binom{y}{1} \leq 0\}$.  Note that 
$\binom01 =0$. Thus $g(1) -1 = 0$ and hence 
$g(1) =1$.\\
\ \\
Thus $(g(1), \ldots, g(5)) = (9,7,3,2,1)$.

Next we apply the unranking proceedure for ${\cal IO}_{4}$ from section 4.3 to 
find $h \in {\cal IO}_{4}$ such that $rank_{{\cal IO}_{4}}(g) = 10$. 
We set $h(4) =0$ and $p_0 =10$.\\
\ \\
Then $10 = 1 \times 3! + 4$ and $4 < 3!$ so that $h(1) = 1$ and $p_1 = 4$.  \\
\ \\
Then $4 = 2 \times 2! + 0$ so that $h(2) = 2$ and $p_2 = 0$.  \\
\ \\
Then $0 = 0 \times 1! + 0$ so that $h(3) = 0$ and $p_3 = 0$.  \\
\ \\
Thus $(h(1),h(2),h(3),h(4)) = (1,2,0,0)$ and $h$ is the direction insertion 
sequence of the permutation $\sigma = 2 \ 4 \ 1 \ 3$.

Finally we apply the unranking proceedure for ${\cal RG}_{7,4}$ from section 
4.4 to find $s = (\alpha_6, \ldots, \alpha_0)$ such that $rank_{{\cal RG}_{7,4}}(s) =100$. \\
\ \\
First we set $\alpha_6  =0$ and $r =100$.\\
\ \\
Then $m_5 = \alpha_6 =0$ and by Table 2, $E^{(4)}(5,0) =65$. Thus the 
maximum value of $s$ such that $s \leq m_5+1$ and $s E^{(4)}(5,0) \leq 100$ is 
$s = 1$. Thus $\alpha_5 =1$ and we set $r = 100 -65 = 35$.\\
\ \\
Then $m_4 = max(\{\alpha_6,\alpha_5\}) =1$ and by Table 2, $E^{(4)}(4,1) =55$. 
Thus the 
maximum value of $s$ such that $s \leq m_4+1$ and $s E^{(4)}(4,1) \leq 35$ is 
$s = 0$. Thus $\alpha_4 =0$ and we set $r = 35$.\\
\ \\
Then $m_3 = max(\{\alpha_6,\alpha_5,\alpha_4\}) =1$ and by Table 2, $E^{(4)}(3,1) =9$. 
Thus the 
maximum value of $s$ such that $s \leq m_3+1$ and $s E^{(4)}(3,1) \leq 35$ is 
$s = 2$. Thus $\alpha_3 =2$ and we set $r = 35 -(2 \cdot 9) = 17$.\\
\ \\
Then $m_2 = max(\{\alpha_6,\ldots,\alpha_3\}) =2$ and by Table 2, $E^{(4)}(2,2) =7$. 
Thus the 
maximum value of $s$ such that $s \leq m_2+1$ and $s E^{(4)}(2,2) \leq 17$ is 
$s = 2$. Thus $\alpha_2 =2$ and we set $r = 17 -(2 \cdot 7) = 3$.\\
\ \\
Then $m_1 = max(\{\alpha_6,\ldots,\alpha_2\}) =2$ and by Table 2, $E^{(4)}(1,2) =1$. 
Thus the 
maximum value of $s$ such that $s \leq m_1+1$ and $s E^{(4)}(1,2) \leq 3$ is 
$s = 3$. Thus $\alpha_1 =3$ and we set $r = 3 -(3 \cdot 1) = 0$.\\
\ \\
Then $m_0 = max(\{\alpha_6,\ldots,\alpha_1\}) =3$ and by Table 2, $E^{(4)}(1,3) =1$. 
Thus the 
maximum value of $s$ such that $s \leq m_0+1$ and $s E^{(4)}(1,2) \leq 0$ is 
$s = 0$. Thus $\alpha_1 =0$.\\
\ \\
Thus $s = (0,1,0,2,2,3,0)$ which corresponds to the set partition 
$\langle (1,3,7), (2), (4,5), (6) \rangle$. 

For Step III, we want to construct $f \in {\cal F}_9$ such 
that 
$$Seq(f) = \langle (9,7,3,2,1),(1,2,0,0),(0,1,0,2,2,3,0)\rangle.$$
Now $\pi^-(f) = \langle (1,3,7), (2), (4,5), (6) \rangle$
Since $\sigma = 2 \ 4 \ 1 \ 3$, the ordered partition 
$\pi(f) = \langle (2), (6), (1,3,7), (4,5) \rangle$ and the 
ordered partition $\pi^+(f) = \langle (3), (7), (2,4,8), (5,6)\rangle$. 
It follows that 
$$Table(f) = \langle\emptyset, \emptyset, \emptyset, (3),(7),(2,4,8), \emptyset,(5,6),\emptyset\rangle.$$
Thus $f$ is specified by the following table.\\
\ \\ 
{\normalsize 
\begin{tabular}{|l|l|l|l|l|l|l|l|l|l|}
\hline 
$f^{-1}(1)$ & $f^{-1}(2)$ & $f^{-1}(3)$ & $f^{-1}(4)$ & $f^{-1}(5)$ & 
$f^{-1}(6)$ & $f^{-1}(7)$ & $f^{-1}(8)$ & $f^{-1}(9)$ \\
\hline
$\emptyset$ & $\emptyset$  & $\emptyset$ & $\{3\}$ & $\{7\}$  & 
 $\{,2,4,8\}$ & $\emptyset$& $\{5,6\}$  & $\emptyset$ \\
\hline 
\end{tabular}}\\
\ \\
The graph of $f$ and the graph of $\Theta(f)$ are pictured in Figure \ref{figure:rank600K}.

\fig{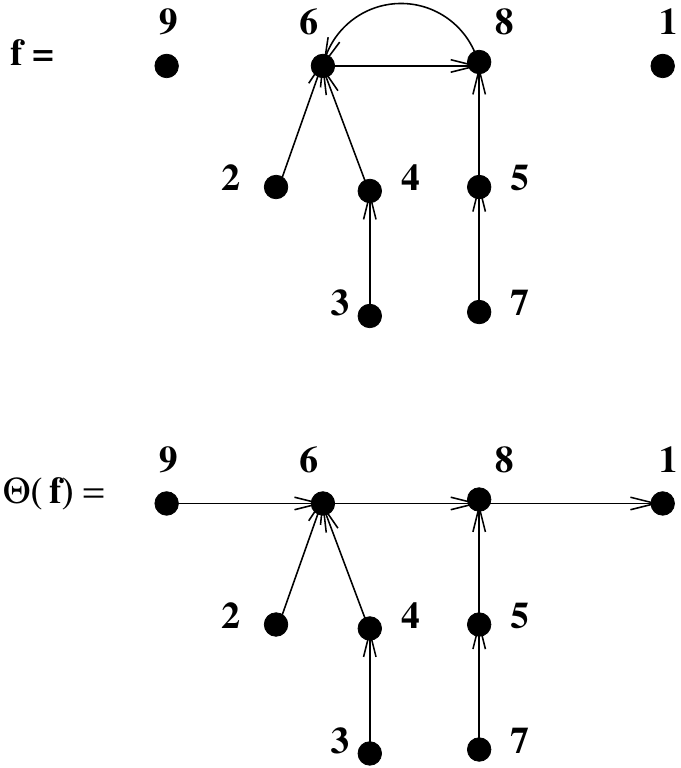}{The element of rank 600,000 in $\vec{C}_{9,1}^5$.}

The ranking and unranking algorithms for $\vec{C}_{n,1}^{k,L}$ where 
$L$ is a subset of $[n]$ of size $k$ is essentially the same as 
the ranking and unranking algorithms for $\vec{C}_{n,1}^k$ except 
that the decreasing function $g \in {\cal DF}_{n,k}$ is precified by $L$. 
That is, if $L = \{l_1 < \ldots < l_k\}$, then 
$g = (g(1), \ldots, g(k)) = (l_k, l_{k-1}, \ldots, l_1)$. 
Now  
$$|\vec{C}_{n,1}^{k,L}| = (n-k)! \times S_{n-2,n-k}.$$
\ \\
{\bf Ranking Algorithm for $\vec{C}_{n,1}^{k,L}$.}\\
\ \\
{\bf Step 1}  Given $T \in \vec{C}_{n,1}^k$, construct $f = \Theta^{-1} (T)$.\\
\ \\
\ \\
{\bf Step 2.}  Find $h \in {\cal IO}_{n-k}$, and 
$s \in {\cal RG}_{n-2,n-k}$ such that 
$$Seq(f) = \langle (g(1), \ldots g(k)), (h(1), \ldots, h(n-k)),(s(1), \ldots, s(n-2))\rangle.$$ 
\ \\
\ \\
{\bf Step 3}.  Use the algorithms of section 4 to compute 
$rank_{{\cal IO}_{n-k}}(h)$  and 
$rank_{{\cal RG}_{n-2,n-k}}(s)$.\\
\ \\
\ \\
{\bf Step 4} It then follows from our analysis in section 4 that 
\begin{equation}\label{Lex1}
rank_{\vec{C}_{n,1}^{k,L}} = 
rank_{{\cal IO}_{n-k}}(h) \cdot S_{n-2,n-k} + 
rank{{\cal RG}_{n-2,n-k}}(s). 
\end{equation}
 Thus if $L = \{3,5,6,8,9\}$, $n =9$, and 
$T$ is the tree pictured in Figure 3, then it follows from our 
previous calculations when we computed $rank_{\vec{C}_{9,1}^5}$, that  
\begin{equation}\label{Lex16}
rank_{\vec{C}_{9,1}^{5,\{3,5,6,8,9\}}}(T) = (7 \times 350) + 46 = 2,496
\end{equation}
Thus the $T$ pictured in Figure 3 is the tree with rank 2,496 in 
$\vec{C}_{9,1}^{5,\{3,5,6,8,9\}}$.

The unranking algorithm for $\vec{C}_{n,1}^{k,L}$ again is essentially 
the same as the unranking algorithm for $\vec{C}_{n,1}^k$ except that 
we do not have to find the decreasing function $g$ since it is prespecified 
by $L$. 
\ \\
\ \\
{\bf Unranking Algorithm for $\vec{C}_{n,1}^k$.}\\
Problem:  Given $r$ such that $0 \leq r < (n-k)! \times S_{n-2,n-k}$, find $T \in \vec{C}_{n,1}^{k,L}$ such that 
$rank_{\vec{C}_{n,1}^{k,L}}(T) =r$.\\
\ \\  
{\bf Step I} Find $r_1$ and $r_2$ such that 
\begin{equation}\label{LEX1}
r =  r_1 \times S_{n-2,n-k} + r_2 \ \mbox{where} \ 0 \leq r_2  < S_{n-2,n-k}.
\end{equation}
\ \\
{\bf Step II} Find $h \in {\cal IO}_{n-k}$ and 
$s \in {\cal RG}_{n-2,n-k}$ such that 
\begin{eqnarray}\label{LEX2}
r_1 &=& rank_{{\cal IO}_{n-k}}(h) \\
r_2 &=& rank_{{\cal RG}_{n-2,n-k}}(s) 
\end{eqnarray}
\ \\
{\bf Step III} Find $f \in {\cal F}_n$ such that 
\begin{equation}\label{LEX3}
Seq(f) = \langle (g(1), \ldots, g(k)), (h(1), \ldots, h(n-k)), (s(1), \ldots, s(n-2))\rangle.
\end{equation}
\ \\
{\bf Step IV} Let $T = \Theta(f)$.
\ \\

For an example of our unranking algorithms for $\vec{C}_{n,1}^{k,L}$, let 
$L = \{3,4,6,8,9\}$ so that $g = (9,8,6,5,3)$. Suppose 
that we want to find the tree $T \in \vec{C}_{9,1}^{5,L}$ whose rank is 
6,000. 
In this case, $(9-5)!S_{7,4} = 8400$ and $S_{7,4} = 350$.  Thus for Step I, 
we find that 
\begin{equation}\label{LEX4}\
6,000=  17 \times 350 + 50.
\end{equation}
Thus $r_1 = 17$ and $r_2 = 50$.

For Step II, first we apply the uranking proceedure for ${\cal IO}_{4}$ from section 4.3 to 
find $h \in {\cal IO}_{4}$ such that $rank_{{\cal IO}_{4}}(g) = 17$. 
We set $h(4) =0$ and $p_0 =17$.\\
\ \\
Then $17 = 2 \times 3! + 5$ and $5 < 3!$ so that $h(1) = 2$ and $p_1 = 5$.  \\
\ \\
Then $5 = 2 \times 2! + 1$ so that $h(2) = 2$ and $p_2 = 1$.  \\
\ \\
Then $1 = 1 \times 1! + 0$ so that $h(3) = 1$ and $p_3 = 0$.  \\
\ \\
Thus $(h(1),h(2),h(3),h(4)) = (2,2,1,0)$ and $h$ is the direction insertion 
sequence of the permutation $\sigma = 3 \ 4 \ 2 \ 1$.

Then we apply the unranking proceedure for ${\cal RG}_{7,4}$ from section 
4.4 to find $s = (\alpha_6, \ldots, \alpha_0)$ such that $rank_{{\cal RG}_{7,4}}(s) =50$. \\
\ \\
First we set $\alpha_6  =0$ and $r =50$.\\
\ \\
Then $m_5 = \alpha_6 =0$ and by Table 2, $E^{(4)}(5,0) =65$. Thus the 
maximum value of $s$ such that $s \leq m_5+1$ and $s E^{(4)}(5,0) \leq 50$ is 
$s = 0$. Thus $\alpha_5 =0$ and we set $r = 50$.\\
\ \\
Then $m_4 = max(\{\alpha_6,\alpha_5\}) =0$ and by Table 2, $E^{(4)}(4,0) =10$. 
Thus the 
maximum value of $s$ such that $s \leq m_4+1$ and $s E^{(4)}(4,1) \leq 50$ is 
$s = 1$. Thus $\alpha_4 =1$ and we set $r = 50-10 =40$.\\
\ \\
Then $m_3 = max(\{\alpha_6,\alpha_5,\alpha_4\}) =1$ and by Table 2, $E^{(4)}(3,1) =9$. 
Thus the 
maximum value of $s$ such that $s \leq m_3+1$ and $s E^{(4)}(3,1) \leq 40$ is 
$s = 2$. Thus $\alpha_3 =2$ and we set $r = 50 -(2 \cdot 9) = 22$.\\
\ \\
Then $m_2 = max(\{\alpha_6,\ldots,\alpha_3\}) =2$ and by Table 2, 
$E^{(4)}(2,2) =7$. 
Thus the 
maximum value of $s$ such that $s \leq m_2+1$ and $s E^{(4)}(2,2) \leq 22$ is 
$s = 3$. Thus $\alpha_2 =3$ and we set $r = 22 -(3 \cdot 7) = 1$.\\
\ \\
Then $m_1 = max(\{\alpha_6,\ldots,\alpha_2\}) =3$ and by Table 2, 
$E^{(4)}(1,3) =4$. 
Thus the 
maximum value of $s$ such that $s \leq m_1+1$ and $s E^{(4)}(1,3) \leq 1$ is 
$s = 0$. Thus $\alpha_1 =0$ and we set $r = 1$.\\
\ \\
Then $m_0 = max(\{\alpha_6,\ldots,\alpha_1\}) =3$ and by Table 2, $E^{(4)}(0,3) =1$. 
Thus the 
maximum value of $s$ such that $s \leq m_0+1$ and $s E^{(4)}(1,2) \leq 1$ is 
$s = 1$. Thus $\alpha_1 =1$.\\
\ \\
Thus in our case, $s = (0,0,1,2,3,0,1)$ which corresponds to the set partition 
$\langle (1,2,6), (3,7), (4), (5) \rangle$. 

For Step III, we want to construct $f \in {\cal F}_9$ such 
that 
$$Seq(f) = \langle (9,8,6,5,3),(2,2,10),(0,0,1,2,3,0,1)\rangle.$$
Now $\pi^-(f) = \langle (1,2,6), (3,7), (4), (5) \rangle$. 
Since $\sigma = 3 \ 4 \ 2 \ 1$, the ordered set partition 
$\pi(f) = \langle (4), (5), (3,7), (1,2,6) \rangle$ and the 
ordered set partition $\pi^+(f) = \langle (5), (6), (4,8), (2,3,7)\rangle$. 
It follows that 
$$Table(f) = \langle (5),(6),\emptyset, (4,8), \emptyset, \emptyset, (2,3,7), \emptyset, \emptyset \rangle.$$
Thus $f$ is specified by the following table.\\
\ \\ 
{\normalsize 
\begin{tabular}{|l|l|l|l|l|l|l|l|l|l|}
\hline 
$f^{-1}(1)$ & $f^{-1}(2)$ & $f^{-1}(3)$ & $f^{-1}(4)$ & $f^{-1}(5)$ & 
$f^{-1}(6)$ & $f^{-1}(7)$ & $f^{-1}(8)$ & $f^{-1}(9)$ \\
\hline
$\{5\}$ & $\{6\} $  & $\emptyset$ & $\{4,8\}$ & $\emptyset$   & 
 $\emptyset$  & $\{2,3,7\}$ & $\emptyset$  & $\emptyset$ \\
\hline 
\end{tabular}}\\
\ \\
The graph of $f$ and the graph of $\Theta(f)$ are pictured in Figure \ref{figure:rank6000}.

\fig{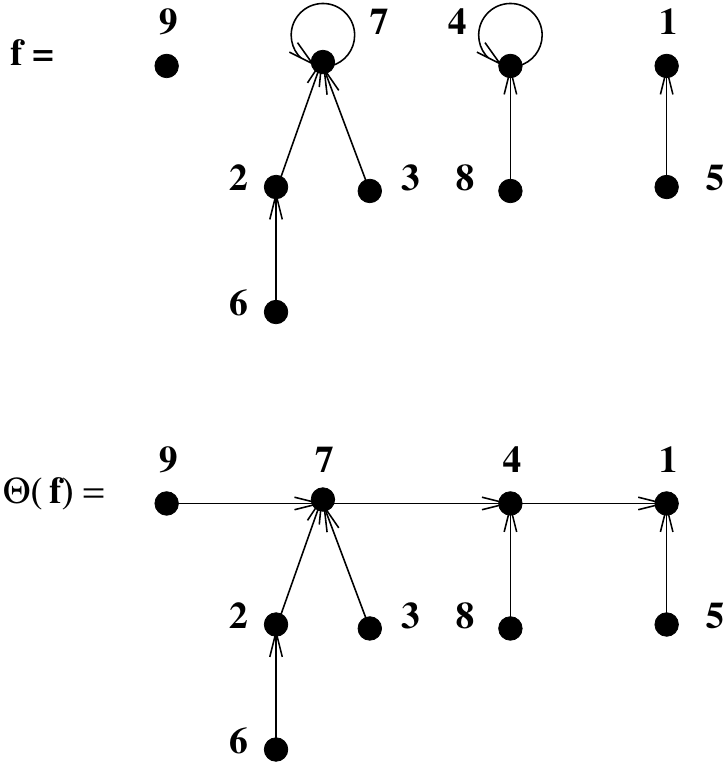}{The element of rank 6,000 in $\vec{C}_{9,1}^{5,\{3,5,6,8,9\}}$.}

Since the ranking and unranking algorithms for 
$\vec{C}_{n,1}^{k,L}$ are essentially a subset of the ranking and unranking 
algorithms for $\vec{C}_{n,1}^k$, it follows that both the ranking and uranking algorithms for $\vec{C}_{n,1}^{k,L}$ require $O(n^2\log(n)$ steps.

\end{document}